\documentclass{amsart}
\usepackage{amsfonts}
\usepackage{amssymb}

\usepackage[all,poly,graph]{xy}

\font\got=eufm10
\font\gots=eufm10 at 7pt
\def\oku{C_{\varphi,\varphi^{-1}}}

\def\ov{\overline}
\def\g2{ \hbox{\got g}_2}
\def\f4{\hbox{\got f}_4}
\def\d4{\hbox{\got d}_4}\def\e6{\hbox{\got e}_6}
\def\fs4{\hbox{\gots f}_4}
\def\F4{\hbox{\got F}_4}
\def\Fs4{\hbox{\gots F}_4}

\def\id{\hbox{\rm id}}

\def\L{{\mathcal
L}}

\def\go{\mathop{\hbox{GO}}}
\def\pgo{\mathop{\hbox{PGO}}}

\def\W{\mathcal W}
\def\V{\mathcal V}

\def\K{{\mathbb K}}
\def\X{\hbox{\got X}}
\def\Xs{\hbox{\gots X}}

\def\C{\mathbb C}
\def\Z{\mathbb Z}
\def\R{\mathbb R}
\def\aut{\mathop{\rm
aut}}
\def\der{\mathop{\rm Der}}

\def\Ad{\mathop{\rm Ad}}
\def\sop{\mathop{\rm Supp}}\def\Out{\mathop{\rm Out}}

\def\Out{\mathop{\rm Out}}
\def\In{\mathop{\rm In}}
\def\Int{\mathop{\rm Int}}
\def\s{\mathop{\rm sl}}

\def\GL{\mathop{\rm GL}}
\def\PGL{\mathop{\rm PGL}}

\def\span#1{\langle #1 \rangle}

\def\sop{\mathop{\hbox{Supp}}}
\def\diag{\mathop{\rm diag}}
\def\a{\alpha}
\def\b{\beta}
\def\al{\ifcase\xypolynode\or F \or A\or B\or C\or D\or G\fi}
\def\ala{\ifcase\xypolynode\or a \or b\or c\or d\or g\or f\fi}

\def\si{\sigma}
\def\To{\hbox{\got T}}
\def\So{\hbox{\got S}}
\def\Co{\hbox{\got C}}
\def\No{\hbox{\got N}}

\def\tor#1{\To^{\langle #1\rangle}}
\def\tors#1{T^{\langle #1\rangle}}
\def\sor#1{\So^{\langle #1\rangle}}
\newtheorem{te}{Theorem}
\newtheorem{pr}{Proposition}
\newtheorem{lm}{Lemma}
\newtheorem{co}{Corollary}
\newtheorem{re}{Remark}
\title[Fine Gradings on $\mathfrak d_4$]{Fine Gradings on the exceptional Lie algebra $\mathfrak d_4$}
\author{Cristina Draper}
\thanks{\ The first and second authors are partially supported by the
MCYT grant MTM2007-60333 and by the JA grants FQM-336, FQM-1215
and FQM-2467. The third author is partially supported by the MCYT
grant MTM2007-60016, and by the JA grants FQM-213 and
P07-FQM-2863.}
\address{Cristina Draper Fontanals: Departamento de
Matemtica Aplicada\\ Campus de El Ejido, S/N, 29013 Mlaga,
Spain\\cdf@uma.es}
\author{C\'{a}ndido Mart\'\i n}
\address{C\'{a}ndido Mart\'\i n Gonzlez:
Departamento de \'Algebra, Geometr\'{\i}a y Topolog\'{\i}a\\
Campus de Teatinos, S/N, 29080 M{\'a}laga,
Spain\\candido@apncs.cie.uma.es}
\author{Antonio Viruel}
\address{Antonio Viruel Arbizar:
Departamento de \'Algebra, Geometr\'\i a y Topolog\'\i a\\ Campus
de Teatinos, S/N, 29080 M{\'a}laga, Spain\\viruel@agt.cie.uma.es}

\begin{document}
\maketitle

  \begin{abstract}
We describe all the fine group gradings, up to equivalence, on the
Lie algebra $\mathfrak d_4$. This problem is equivalent to finding the maximal abelian
diagonalizable subgroups of the automorphism group of $\mathfrak
d_4$.
We prove that there are fourteen by using
two different viewpoints. The first approach is computational:
we get a full description of the gradings
by using a particular implementation of the automorphism group
of the Dynkin diagram of $\mathfrak d_4$ and some
algebraic groups stuff. The second approach, more qualitative,
emphasizes some algebraic aspects, as triality, and it is mostly
devoted to gradings involving the   outer automorphisms of
order three.
\end{abstract}

\section{Introduction}

The increasing mathematical activity around gradings on Lie and
other kind of algebras is a phenomenon which is in the background
of the choices of bases and of maximal sets of quantum observables
with additive quantum numbers, as pointed out by Patera et al.  in
a series of papers on the subject \cite{Zass, pat1,  Phys3, o4}.
The Cartan decomposition of a semisimple Lie algebra is one of the
most interesting fine gradings (in fact the unique toral one), with
a heavy influence in the structure of Lie algebras and in
representation theory. But this influence invades also the nearby
fields of particle physics via the usual identification of
observables with generators in a Cartan subalgebra, and particles
living comfortably in the root spaces of a suitable Lie algebra.
So, for instance, it is possible to model the strong interaction of
nature by means of the $\g2$ exceptional Lie algebra. The
possibility of using this algebra for describing hypercharge and
isospin third component, for a series of $14$ elementary particles
(quarks $u$, $d$, $s$ and mesons $\pi^+$, $K^+$, $K_0$, $\pi^0$
together with their antiparticles), was highlighted by   Gunaydin
and   Gursey in the early seventies (see \cite{new1}.)  In this
nice description, two generators of the Cartan subalgebra are
identified with the observables of hypercharge and isospin third
component. These act simultaneously ad-diagonally on the root
spaces (so that the roots give the quantum numbers.)  And root
space generators are the elementary particle mentioned above. So
this is a physical picture of the fine toral grading of $\g2$. In
a more recent development, this scheme is repeated when describing
the strong interaction force by means also of $\g2$ (see for
instance \cite[p.~5]{new2}.)  Thus, strong interaction in nature
may be described as the fine toral grading of $\g2$ with six
gluons as long roots ($g^{r\bar g}$, $g^{\bar rg}$, $g^{r\bar b}$,
$g^{\bar rb}$, $g^{\bar gb}$ and $g^{g\bar b}$), and color quarks
$q^r$, $\bar q^r$, $q^g$, $\bar q^g$, $q^b$ and $\bar q^b$ as
short roots, so that each antiquark is the opossite root of the
given quark. Along the lines used to describe strong interactions,
other Lie algebras are being used to tentatively describe the rest
of the forces in nature. In many of these attempts,  fine
gradings appear as a common feature.

Other independent motivations toward gradings come from the
theory of contraction of Lie algebras. As claimed in \cite{new3},
contractions are important in physics because they explain in
terms of Lie algebras why some theories arise as a limit of more
exact  theories.  For example, the passage from the Poincar\'e
algebra to the Galilei algebra (as $c\to\infty$.)  Contraction
consists in multiplying the generators of the symmetry by
contraction parameters,  such that when these parameters reach
some singularity point one obtains a non-isomorphic Lie algebra
with the same dimension. Graded contractions are a key ingredient
when studying contractions that keep ``undeformed" certain
subalgebra (e.g. see \cite{cont}.)

Gradings on Lie algebras are also of interest in the setting of   Jordan algebras: given a group
grading on $L=\sum_{g\in G}L_g$  and a nontrivial element $g\in G$, any element
$k\in L_{g^{-1}}$ gives rise to a Jordan algebra structure in $L_g$ by
defining $x\circ y=(xk)y$. This also provides a bridge between group gradings and physics,
by the classical work of Jordan \cite{jordan}.

The Lie algebra  of type $\d4$, realized as the algebra of
skew-symmetric matrices $\frak o(8,\K)$, is considered as
exceptional by many authors. In spite of being a member of the
family of classical Lie algebras $\hbox{\got d}_n$, it is the only
one enjoying the benefits of a triality automorphism coming from
the special triangular symmetry of its Dinkyn diagram. It can be
said that this is the most symmetric algebra in the family.  But
more symmetries imply more gradings.  This may be the reason why
the study of gradings on $\d4$ does not match the general scheme
of gradings in the rest of the algebras $\hbox{\got d}_n$. For
instance, the gradings on the algebras $\mathfrak o(n,\K)$ are
computed in \cite{Ivan2} and \cite{LGII}, but only for $n\ne8$.

The main notions about Lie gradings are given in \cite{Zass} by
Patera and Zassenhaus, with a certain  rectification in
\cite{Alb2} about the existence on a grading group, which must be
assumed. Continuation of that work are the papers \cite{LGII},
which deals with the gradings on the classical Lie algebras of
types $\mathfrak{sl}(n,\K)$ and ${\mathfrak o}(n,\K)$, and
\cite{LGIII}, where the real case is considered. An alternative
line of working is followed by Shestakov and Bahturin in
\cite{Ivan2}, but again by using tools of associative algebras.
The first work containing a treatment  of the gradings on a
exceptional Lie algebra is \cite{g2}, which describes the gradings
on $\mathfrak{g}_2$ by taking the octonions as a starting point
(see also \cite{otrog2}.)  But the techniques used in that case are
not enough to obtain the gradings on $\f4$, so that the same
authors develop some computational techniques in \cite{f4} to
obtain a complete description of the nontoral gradings on $\f4$.
Such tools are applied successfully in this paper to obtain the
fine gradings on $\d4$, but we want to remark the existence of
alternative methods which could be useful possibly in algebras of
bigger rank.\smallskip

The paper is organized as follows. In Section 2 we give a quick
review of the main notions relative to Lie gradings and their
translation in terms of groups. Section 3 deals with the problem
of finding fine gradings of $\d4$ from a computational approach.
Since every MAD-group of $\aut\d4$ lives in the normalizer of a
maximal torus, the first task is to fix a maximal torus. Now, to
compute the normalizer we need one automorphism extending each
element in the isometry group of the root system. Next, in
Theorem~\ref{losmads}, we give an explicit expression of all the
MAD's in terms of these elements, which, in particular, gives us
their grading groups and the types of the gradings. Notice that
the MAD's of $\aut\d4$ that map into a $2$-sylow of the group of
components of $\aut \d4$, denoted by $\Out\d4\cong S_3$, are
essentially described in \cite{LGII} (at least, how to get them,
although without specifying the properties, description or grading
groups.)  That is why, in Theorem~\ref{losmads}, we prove that
there are exactly three MAD's which map onto the $3$-sylow of
$\Out\d4$, and we add a description of the remaining MAD's for its
possible applications (a complete matrix description of all the
fine group gradings on $\d4$ is given in \cite{acta}.)

In Section~4 we obtain the same results from a completely
different approach, that does not require any computer
calculation. This allows to have a quick idea of which are these
three MAD's containing outer $3$-automorphisms. In fact, we can
understand them deeply (with the concrete homogeneous components
of the induced grading) only by knowing the gradings on the
simpler Lie algebras $\hbox{\got g}_2$ and $\hbox{\got a}_2$ (or
easier still, the gradings on the octonion algebra and on
$M_{3\times3}(\C)$.) This suggests an inductive process based in
the knowledge the gradings on Lie algebras of less rank to compute
the unknown gradings on the Lie algebras $\hbox{\got e}_6$,
$\hbox{\got e}_7$ and $\hbox{\got e}_8$ (applied to $\hbox{\got
e}_6$ with success in \cite{casiya}.)

\section{Preliminary notions}

  If $V$ is a finite-dimensional Lie  algebra and $G$ is an
abelian group, we shall say that the decomposition $V=\oplus_{g\in
G}V_g$ is a \emph{$G$-grading} whenever for all $g,h\in G$,
$V_gV_h\subset V_{gh}$ and $G$ is generated by the set
$\sop(G):=\{g\in G\colon V_g\ne0\}$, called the \emph{support} of
the grading. We say that two gradings $V=\oplus_{g\in
G}X_g=\oplus_{g'\in G'}Y_{g'}$ are \emph{equivalent} if the sets
of homogeneous subspaces are the same up to isomorphism, that is,
there are an automorphism $f\in\aut(V)$ and a bijection between
the supports $\alpha\colon \sop(G)\to \sop(G')$ such that
$f(X_g)=Y_{\a(g)}$ for any $g\in \sop(G)$.  A convenient invariant
for equivalence is that of \emph{ type}. Suppose we have a grading
on a finite dimensional algebra, then for each positive integer
$i$ we will denote (following \cite{Hesse}) by $h_i$ the number of
homogeneous components of dimension $i$. In this case we shall say
that the grading is {of type} $(h_1,h_2,\dots,h_l )$, for $l$ the
greatest index such that $h_l\ne0$.   Of course, the number $\sum_i
i h_i$ agrees with the dimension of the algebra. We shall say that
the $G$-grading is a \emph{refinement} of the $G'$-grading if and
only if each homogeneous component $Y_{g'}$ with $g'\in G'$ is a
direct sum of some homogeneous components $X_g$. A   grading is
{\em fine} if its unique refinement is the given grading. Our
objective is to classify fine gradings up to equivalence.

The ground field $\K$ will be supposed to be algebraically closed
and of characteristic zero  throughout this work. Notice that the
group of automorphisms of the algebra $V$ is an algebraic linear
group.  There is a deep relationship between gradings on $V$ and
quasitori of the group of automorphisms $\aut (V)$, according to
\cite[\S 3, p.\,104]{enci}.
% If $S$ is a finitely generated
%abelian group, then its group of characters $\X(S)=
%\hom(S,k^\times)$ is a quasitorus and reciprocally, the group of
%characters of a quasitorus turns out to be a finitely generated
%abelian group.
If $V=\oplus_{g\in G}V_g$ is a $G$-grading, the map
$\psi\colon\X(G)=\hom(G,\K^\times)\to\aut(V)$ mapping each
$\alpha\in\X(G)$ to the automorphism $\psi_\alpha\colon V\to V$
given by $V_g\ni x\mapsto \psi_\alpha(x):=\alpha(g)x$ is a group
homomorphism. Since $G$ is finitely generated, then $\psi(\X(G))$
is a quasitorus. And conversely, if $Q$ is a quasitorus and
$\psi\colon Q\to\aut(V)$ is a homomorphism, $\psi(Q)$ is  formed
by semisimple automorphisms
 and we have a $\X(Q)$-grading
$V=\oplus_{g\in{\Xs}(Q)}V_g$ given by $ V_g=\{x\in V\colon
\psi(q)(x)=g(q)x\ \forall q\in Q\} $, with $\X(Q)$  a finitely
generated abelian group. If $V=\oplus_{g\in G}V_g$ is a
$G$-grading,   the set of automorphisms of $V$ such that every
$V_g$ is contained in some eigenspace is an abelian group formed
by semisimple automorphisms, which contains to $\psi(\X(G))$.
%y el grupo de caracteres de dicho grupo es el grupo universal de la graduacion,
% o sea el maxiaml (respecto a la relacion de contenido) entre los que
%producen graduaciones equivalentes a una dada
The grading is fine if and only if such set is a maximal abelian
subgroup of semisimple elements, usually called a \emph{ MAD
}(\lq\lq maximal abelian diagonalizable")-group. It is convenient to
observe that the number of conjugacy classes of MAD-groups of
$\aut(V)$ agrees with the number of equivalence classes of fine
gradings on $V$.

\section{Computational approach}

 First of all we must invoke a version of
the  Borel-Serre  theorem (\cite[Theorem~3.15, p.\,92]{Platonov})
asserting that a
  supersolvable subgroup
  of semisimple elements in an algebraic group is contained
   in the normalizer of some maximal torus. In particular, this
   can be applied to finitely generated abelian groups.
   As we are able to implement a concrete maximal torus and its
   normalizer in a computer, we will use computational methods to find  the
   maximal quasitori, which, up to conjugation,
    live in that normalizer.

Fix the symmetric matrix $C=\begin{pmatrix}
0&I_4\\I_4&0\end{pmatrix}$ and consider the Lie algebra
$\L:=\{x\in\text{Mat}_{8\times 8}(\K)\colon x^tC=-Cx\}$,
 of type $\d4$. Take the Cartan subalgebra $\mathfrak h$ formed
by the diagonal matrices of $\L$. Let $\L=\sum_{\alpha\in\mathfrak
h^*} L_\alpha$ be the root decomposition relative to $\mathfrak
h$, that is, $L_\alpha=\{x\in\L\colon [h,x]=\alpha(h)x\quad\forall
h\in\mathfrak h\}$, and $\Phi=\{\alpha\in\mathfrak h^*-\{0\}\colon
L_\alpha\ne0\}$  is the root system. Let $e_{i,j}$ denote the
elementary matrix whose all entries are trivial but the
$(i,j)$-entry, which is $1$. If we write
$h_i=e_{i,i}-e_{i+4,i+4}$, then $\{h_i\colon i=1,\dots, 4\}$ is a
basis of $\mathfrak h$.

If $h=\sum_{i=1}^4 w_ih_i$ is an arbitrary element in $\mathfrak
h$, and we define $\alpha_i\colon\mathfrak h\to \K$  by
$\alpha_i(h)=w_i-w_{i+1}$, for $i=1,2,3$, and
$\alpha_4(h)=w_3+w_4$, then
$\Delta=\{\alpha_1,\alpha_2,\alpha_3,\alpha_4\}$ is a basis of
  $\Phi$. Indeed, if $b_{i,j}:=e_{j,i}-e_{i+4,j+4}$,
$c_{i,j}:=e_{j,i+4}-e_{i,j+4}$ and $d_{i,j}:=e_{i+4,j}-e_{j+4,i}$,
we have $[h,b_{i,j}]=(w_j-w_i)b_{i,j}$, $
[h,c_{i,j}]=(w_j+w_i)c_{i,j}$ and $[h,d_{i,j}]=(-w_i-w_j)d_{i,j}$.
Thus we can choose $\mathcal{B}=\{h_i\colon
i=1,\dots,4\}\cup\{b_{i,j}\colon i\ne
j,\,i,j=1,\dots,4\}\cup\{c_{i,j},d_{i,j}\colon
i<j,\,i,j=1,\dots,4\}$ a basis formed by root vectors, with $b$'s,
$c$'s and $d$'s ordered following the rows in the next array
\begin{equation}\label{base}
\begin{array}{lll}
b_{2,1}\in L_{\alpha_1}&b_{3,2}\in L_{\alpha_2}&b_{4,3}\in
L_{\alpha_3}\\
c_{3,4}\in L_{\alpha_4}&b_{3,1}\in
L_{\alpha_1+\alpha_2}&b_{4,1}\in L_{\alpha_1+\alpha_2+\alpha_3}\\
b_{4,2}\in L_{\alpha_2+\alpha_3}&c_{4,1}\in
L_{\alpha_1+\alpha_2+\alpha_4}&c_{4,2}\in L_{\alpha_2+\alpha_4}\\
c_{2,1}\in L_{\alpha_1+2\alpha_2+\alpha_3+\alpha_4}&c_{3,1}\in
L_{\alpha_1+\alpha_2+\alpha_3+\alpha_4}&c_{2,3}\in
L_{\alpha_2+\alpha_3+\alpha_4}
\end{array}
\end{equation}
followed by the opposite roots (take $b_{j,i}\in L_{-\alpha}$ when
$b_{i,j}\in L_\alpha$, and $d_{ i,j}\in L_{-\alpha}$ when
$c_{i,j}\in L_\alpha$) in the same order.

Take $\To:=\{t\in\aut\L\colon t\vert_{\mathfrak h}=\id\}$. This is
a maximal torus of $\aut\L$ such that each element acts diagonally
on the root spaces. More precisely, if the automorphism acts with
eigenvalues $\alpha,\beta,\gamma$ and $ \delta$ in $L_{\alpha_i}$
respectively for $i=1,\dots,4$, then its matrix relative to the
basis $\mathcal{B}$ is
$$
\begin{array}{c}
\text{diag}\{1,1,1,1,\alpha,\beta,\gamma,\delta,\alpha\beta,\alpha\beta\gamma,
\beta\gamma,
\alpha\beta\delta,\beta\delta,\alpha\beta^2\gamma\delta,\alpha\beta\gamma\delta,\beta\gamma\delta,\\
\frac1\alpha,\frac1\beta,\frac1\gamma,\frac1\delta,\frac1{\alpha\beta},\frac1{\alpha\beta\gamma},\frac1{\beta\gamma},
\frac1{\alpha\beta\delta},\frac1{\beta\delta},\frac1{\alpha\beta^2\gamma\delta},
\frac1{\alpha\beta\gamma\delta},\frac1{\beta\gamma\delta} \}
\end{array}
$$
and the automorphism will be denoted by
$t_{\alpha,\beta,\gamma,\delta}$.

In order to get the normalizer of $\To$, we need to describe the
abstract Weyl group of $\d4$. The
 Cartan matrix of $\d4$ is
$$  \small\begin{pmatrix} 2 & -1 & 0 & 0\cr -1 & 2 &
-1 & -1\cr 0 & -1 & 2 & 0\cr 0 & -1 & 0 & 2
\end{pmatrix}.
$$
So, for
 $E=\sum_{i=1}^4\R\a_i$ the euclidean space
  with the inner product $\langle\ ,\ \rangle$,
 the Weyl group of $\d4$ is the subgroup $\W$ of $\GL(E)$
generated by the (simple) reflections $s_i$ with $i=1,2,3,4$,
given by $s_i(x):=x-\span{x,\a_i}\a_i$, where the Cartan integers
$\langle\cdot,\cdot\rangle$ are extracted as usual from the Cartan
matrix above. Identifying $\GL(E)$ to $\GL(4,\R)$ by means of the
matrices relative to the $\R$-basis $\Delta$, the reflections
$s_i$ are represented by
$$
\small{\begin{array}{ll} s_1=\begin{pmatrix}-1 & 0 & 0 & 0\cr 1 &
1 & 0 & 0\cr 0 & 0 & 1 & 0\cr 0 & 0 & 0 & 1\end{pmatrix},\
&s_2=\begin{pmatrix}1 & 1 & 0 & 0\cr 0 & -1 & 0 & 0\cr 0 & 1 & 1 &
0\cr 0 & 1 & 0 & 1\end{pmatrix},\\\,\\s_3=\begin{pmatrix}1 & 0 & 0
& 0\cr 0 & 1 & 1 & 0\cr 0 & 0 & -1 & 0\cr 0 & 0 & 0 &
1\end{pmatrix},\ &s_4=\begin{pmatrix}1 & 0 & 0 & 0\cr 0 & 1 & 0 &
1\cr 0 & 0 & 1 & 0\cr 0 & 0 & 0 & -1\end{pmatrix}.\end{array}}
$$
 Now the group of isometries of $\Phi$ is the
semidirect product of the Weyl group with the group of
automorphisms of the Dynkin diagram. These automorphisms come from
permutations $\sigma$ of   $\{1,\dots, 4\}$ such that
$\span{\alpha_i,\alpha_j}=\span{\alpha_{\sigma(i)},\alpha_{\sigma(j)}}$.
Hence, denoting by
$$
\small{ s_5=\begin{pmatrix}0 & 0 & 1 & 0\cr 0 & 1 & 0 & 0\cr 0 & 0
& 0 & 1\cr 1 & 0 & 0 & 0\end{pmatrix},\qquad  s_6=\begin{pmatrix}1
& 0 & 0 & 0\cr 0 & 1 & 0 & 0\cr 0 & 0 & 0 & 1\cr 0 & 0 & 1 &
0\end{pmatrix},}
$$
respectively related to the pictures
%dibujos
\vskip .3cm \hskip 3cm
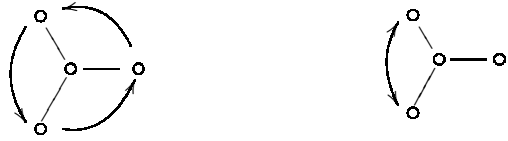
\vbox{
\xy <1cm, 0cm>:
\POS(0,0)*+\xycircle(.07,.07){-}
\ar @{-} +(.33,-.57)
\POS(.4,-.7)*+\xycircle(.07,.07){-}
\ar @{-} +(.65,0)
\POS(1.3,-.7)*+\xycircle(.07,.07){-}
\POS(.35,-.82)\ar @{-} +(-.33,-.57)
\POS(0,-1.5)*+\xycircle(.07,.07){-};
(1.2,-.4);(.3,0,1)
**\crv{(0.9,0.23)}?>*\dir{>} ;
 (-0.2,-1.4);(-0.2,-0.15)
**\crv {(-.6,-.775)}?<*\dir{<};
(0.3,-1.5);(1.25,-.85)
**\crv{(0.9,-1.6)}?>*\dir{>}
\endxy}
\vskip -1.9cm \hskip 8cm
\vbox{
\xy <1cm, 0cm>
\POS(0,0)*+\xycircle(.7,.7){-}
\ar @{-} +(3.5,-5.9)*+
\xycircle(.7,.7){-}
\POS(5,-5.9)
\ar @{-} +(6.5,0)*+\xycircle(.7,.7){-}
\POS(3,-7)
\ar @{-} +(-2.7,-5)
\POS(0,-13)*+
\xycircle(.7,.7){-};
(-2,-1);(-2,-12)
**\crv{(-5,-6)}?>*\dir{>}?<*\dir{<}
\endxy}
\vskip .5cm

\noindent the group of automorphisms of the Dynkin diagram of
$\d4$ is $\langle s_5,s_6\rangle\cong S_3$ and the group
$\aut\Phi=\span{s_i\colon i=1,\dots
6}\equiv\mathcal{V}\cong\W\rtimes S_3$.

 We shall consider $\mathcal{V}\subset \GL(4,\R)$ ordered
lexicographically, that is, first for any two different couples
$(i,j), (k,l)$ such that $i,j,k,l\in\{1,2,3,4\}$ we define
$(i,j)<(k,l)$ if and only if either $i<k$ or $i=k$ and $j<l$, and
second, for any two different matrices $\si=(\si_{ij})$,
$\si'=(\si'_{ij})$ in $\mathcal{V}$,  $\si<\si'$ if and only if
$\si_{ij}<\si'_{ij}$ where $(i,j)$ is the least element (with the
previous order in the couples)  such that $\si_{ij}\ne\si'_{ij}$.
One possible way to compute this group with this particular
enumeration is provided by the following code implemented with
{\sl Mathematica}:
\smallskip

{\parindent=3cm\tt V=Table[$s_i$,\{i,6\}];

a[L\_,x\_]:=Union[L, Table[L[[i]].x,\{i,Length[L]\}],

\hskip 4cm Table[x.L[[i]],\{i,Length[L]\}]]

Do[V=a[V,$s_i$],\{i,6\}]\hskip 1cm \textrm{(3 times repeated)} }
\smallskip

We get a list of $1152=2^73^2$ elements in the table {\tt V} which
is nothing but the   group $\aut\Phi$. We are denoting by
$\sigma_i$ the $i$-th element of $\mathcal{V}$ lexicographically
ordered.
%The following result comes from a straightforward computation
%which may be done with any matrix multiplication software.
Now, for each $\sigma\in\aut\Phi$, it is possible to choose an
automorphism  $\tilde\sigma\in\aut\L$ mapping $L_\alpha$ into
$L_{\sigma(\alpha)}$ such that $\tilde\sigma\vert_{\mathfrak h}$
agrees with $\sigma\in\text{End}(\mathfrak h^*)$ by means of the
identification between $\mathfrak h$ and $\mathfrak h^*$  given by
the Killing form ($h\in\mathfrak h\mapsto K(h,-)\in\mathfrak
h^*$.) Concretely, if $v_\alpha$ denotes the root vector specified
in (\ref{base}), we take as $\tilde \sigma$   the only
automorphism such that $\tilde\sigma(v_\alpha)=v_{\sigma(\alpha)}$
according to the isomorphism theorems in \cite{Humphreysalg} (any
other choice would have been in the form $\tilde\sigma t$ for some
$t\in\To$.) Thus we have a precise description of the normalizer
of the maximal torus as
$$
 N_{\aut\d4}(\To)=\{\tilde\sigma_it_{x,y,z,u}\colon
i=1,\dots,1152,\,x,y,z,u\in \K^\times\}=:\No.
$$
 The point is  that all
the MAD-groups of $\aut\d4$ live in $\No$ (up to conjugation), so
we will be able to find a concrete description of them in terms of
$\tilde\sigma_i$ and $t_{x,y,z,u}$. This approach allows us to
know in detail the homogenous components.

\begin{te}\label{losmads} There are fourteen maximal quasitori in $\aut\d4$.
They are isomorphic to
$$\begin{array}{l}
\Z_2^{n}\qquad
(n=5,6,7)\\ \K^\times\times\Z_2^{n}
\qquad(n=3,4,5)\\(\K^\times)^2\times \Z_2^{n}\qquad (n=2,3)\\
(\K^\times)^3\times\Z_2\\(\K^\times)^4\\\Z_4\times\Z_2^3\\
\Z_3^3,\,\Z_3\times\Z_2^3,\,\Z_3\times(\K^\times)^2\end{array} $$
 and their precise descriptions, jointly with the types of the gradings
induced by them, are summarized in the following table.
\end{te}

\begin{center}
\begin{tabular}{|c|c|c|c|c|}
\hline   & Grading  & Automorphisms generating the group &
Type&$\dim L_e$ \cr &group&&&\cr\hline   $Q_1$&$\Z_2^7$
&$t_{-1,1,1,1},t_{1,1,-1,-1},t_{1,-1,1,1}$ & $(28)$&$0$\cr   &
&$\tilde{\sigma}_{1},\tilde{\sigma}_{3},\tilde{\sigma}_{19},\tilde\sigma_{259}$&&\cr
\hline$Q_2$& $\Z_2^5\times\Z$&
$t_{-1,1,1,1},t_{1,1,-1,-1},t_{1/u,u,1,1}$& $(28)$&$1$\cr &
&$\tilde\sigma_{1},\tilde\sigma_{3},\tilde\sigma_{19}$ &&\cr
\hline $Q_3$& $\Z_2^3\times\Z^2$&
$t_{-1,1,1,1},t_{1/u,u,1,1},t_{1/v,1,v,v}$& $(26,1)$&$2$\cr &
&$\tilde\sigma_{1},\tilde\sigma_{3}$ &&\cr \hline $Q_4$&
$\Z_2\times\Z^3$&
$t_{1/u,u,1,1},t_{1/v,1,v,v},t_{w,1,1,1/w^2},\tilde\sigma_{3}$&
$(25,0,1)$ &$3$\cr \hline $Q_5$& $\Z_2^3\times\Z$&
$t_{-1,1,1,1},t_{1/v,1,v,v},t_{1,-1,1,-1},\tilde\sigma_{49}t_{1,-1,1,1}$&
$(25,0,1)$ &$1$\cr \hline $Q_6$& $\Z_2^3\times\Z_4$&
$t_{1,-1,1,1},t_{-1,-1,1,1},
\tilde\sigma_{259}t_{-1,-1,-1,-1},\tilde\sigma_{7}$& $(24,2)$
&$0$\cr \hline $Q_{7}$& $\Z_2^4\times\Z$&
$t_{-1,1,-1,1},t_{1,-1,1,1},t_{1,1,1,u},\tilde\sigma_{280}
,\tilde\sigma_{634}$& $(28)$ &$1$\cr \hline $Q_{8}$& $\Z_2^6$&
$t_{-1,1,1,1},t_{1,-1,1,1},t_{1,1,-1,-1}$& $(28)$ &$0$\cr &
&$\tilde\sigma_{1}t_{i,1,i,-i},\tilde\sigma_{259}t_{1,1,1,-1},
\tilde\sigma_{243}t_{-1,-i,1,1}$ &&\cr \hline $Q_{9}$& $\Z_2^5$&
$t_{-1,1,1,1},t_{1,-1,1,1},t_{1,1,-1,1},t_{1,1,1,-1},\tilde\sigma_{259}$&
$(24,0,0,1)$ &$0$\cr \hline $Q_{10}$& $\Z_2^2\times\Z^2$&
$t_{-1,1,1,1},t_{1/u,u,1,1},t_{1/v,1,v,v},\tilde\sigma_{1}$&
$(20,4)$ &$2$\cr \hline $Q_{11}$& $ \Z^4$&
$t_{u,1,1,1},t_{1,v,1,1},t_{1,1,w,1},t_{1,1,1,z}$& $(24,0,0,1)$
&$4$\cr \hline $Q_{12}$& $\Z_2^3\times\Z_3$&
$t_{1,-1,1,1},t_{-1,-1,-1,-1},\tilde\sigma_{20}$& $(14,7)$&$0$\cr
\hline $Q_{13}$&$\Z_3\times\Z^2$
&$t_{1,y,1/y^2,1},t_{x,1,1/x^3,x},\tilde\sigma_4$ &
$(26,1)$&$2$\cr \hline
$Q_{14}$&$\Z_3^3$&$t_{1,\omega,\omega,\omega},t_{\omega,1,\omega,1},\tilde\sigma_{59}$
& $(24,2)$ &$0$\cr \hline
\end{tabular}
\end{center}

\begin{center}
\textsc{Table 1}
\end{center}\smallskip

The   zero-component $L_e$ in a fine grading of a Lie algebra $L$
is an abelian subalgebra whose dimension is the dimension of the
quasitorus producing such grading, as showed in
\cite[Prop.\,10]{f4}.
\smallskip

Observe that the type is not enough to determine the isomorphy
classes of fine gradings on $\d4$, in contrast to the gradings on
$\f4$. For instance there are 4 fine gradings with every
homogeneous component one-dimensional.

\begin{re} \rm The correspondence with the MAD-groups which
could be obtained following the lines in \cite{LGII} is given by:
$T_{0,8}^{(0)}\cong Q_{1}$ (conjugated groups),
$T_{2,6}^{(0)}\cong Q_{2}$, $T_{4,4}^{(0)}\cong Q_{3}$,
$T_{6,2}^{(0)}\cong Q_{4}$, $ T_{2,2}^{(1)}\cong Q_{5}$,
$T_{0,2}^{(2)}\cong Q_{6}$, $T_{2,0}^{(2)}\cong Q_{7}$, $
T_{0,1}^{(3)}\cong Q_{8}$, $T_{0,4}^{(1)}\cong Q_{9}$,
$T_{4,0}^{(1)}\cong Q_{10}$, $ T_{8,0}^{(0)}\cong Q_{11}$.
Although that paper does not consider the outer automorphisms of
order 3, it clearly shows how to find explicit expressions for
the generators of the remaining MAD's, since
$\Int\mathfrak{o}_K(8,\C)\cdot\Z_2=\Ad \rm{O}_K(8,\C)$ (for $\Ad
P(x)=PxP^{-1}$.) Each of the $T_{2k,8/m-2k}^m$ is expressed in
terms of a convenient symmetric matrix $K$.\end{re}

Before proving Theorem~\ref{losmads}, we are going to introduce
some tools. Denote by $\pi\colon \No\to\mathcal V$ the group
epimorphism given by $\pi(\tilde\sigma_it_{x,y,z,u})=\sigma_i$.
Notice that the action $\V\times\To\to\To$ given by $\sigma\cdot
t=\tilde\sigma t\tilde\sigma^{-1}\in\To$ does not depend on our
choice of the concrete extension $\tilde\sigma$, but $\si\cdot
t_{x,y,z,u}=t_{x',y',z',u'}$ for
 $$
\begin{matrix}
 x'= & x^{b_{11}}y^{b_{12}}z^{b_{13}}u^{b_{14}},\cr
 y'= & x^{b_{21}}y^{b_{22}}z^{b_{23}}u^{b_{24}},\cr
 z'= & x^{b_{31}}y^{b_{32}}z^{b_{33}}u^{b_{34}},\cr
 u'= & x^{b_{41}}y^{b_{42}}z^{b_{43}}u^{b_{44}},
 \end{matrix}
$$
if $\sigma=(b_{ij})_{i,j=1,\dots,4}\in\V$. Now consider the
quasitori
$$\begin{array}{l}
\tor{j}:=\{t\in\To\colon\tilde\sigma_j \cdot t=t\},\\
Q(j,t_0):=\span{{\tilde\sigma_jt_0}\cup\tor{j} },
\end{array}$$
 if $j\in\{1,\dots
1152\}$ and $t_0\in\To$. Notice that the   maximal quasitori not
contained in $ \Int(\d4)\cdot\Z_2$, that is, the ones related to
the triality automorphism, will be proved to be  precisely
$Q_{12}=Q(20,\id)$, $Q_{13}=Q(4,\id)$ and $Q_{14}=Q(59,\id)$. This
explains the relevance of the considered quasitori. In some
groups, like $\aut\f4$, every quasitorus is a subgroup of some
$Q(j,\id)$ (see \cite{f4}.) Obviously, this is not our case (some of the $Q_i$'s
have more than 5 generators), but anyway,
 any $Q$ minimal nontoral quasitorus (that is, $Q$ is non
toral but it does not contain properly   any nontoral quasitorus
of $\aut\d4$) is always conjugated to a subgroup of some $Q(j,t)$. As
$Q(j,t)$ is isomorphic to some $Q(i,t')$ if $\sigma_j$ is
conjugated to $\sigma_i$, we only have  to consider the indices of
some representatives of the   orbits up to conjugation in $\V$.
There are 139 order 2 elements distributed in 7 orbits, 80 order 3
elements distributed in 3 orbits, 228 order 4 elements in 5
orbits, 464 order 6 elements in 7 orbits, 144 order 8 elements in
only one orbit, and 96 order 12 elements in the same orbit. One
choice of representatives of these orbits, jointly with the
relevant quasitori related to them,  is the following:

\begin{center}
\begin{tabular}{|c|c|c|c|}
\hline order &   representative& $\tor{j}$ & $Q(j,\id)$\cr & of
the orbit&&is  toral?\cr \hline
 1 & 894 & $\To\cong\K^{\times4}$& yes\cr
 2 & 1& $\{t_{\frac{x}{yz},y,z,z}\colon x^2=1\}\cong \K^{\times2}\times\Z_2$& no \cr
 2 & 3& $\{t_{x,y,z,\frac1{x^2y^2z}}\colon x,y,z\in \K^\times\}\cong \K^{\times3}$&  no\cr
 2 & 9& $\{t_{x,y,\frac1{x^2y^2},1}\colon x,y\in \K^\times\}\cong \K^{\times2}$&  no\cr
 2 & 19& $\{t_{x,\frac{y}x,u,u}\colon y^2=u^2=1\}\cong \K^{\times}\times\Z_2^2$&  no\cr
 2 & 49& $\{t_{\frac{x}z,y,z,yz}\colon x^2=y^2=1\}\cong \K^{\times}\times\Z_2^2$&  no\cr
 2 & 259& $\{t_{x,y,z,u}\colon x^2=y^2=z^2=u^2=1\}\cong  \Z_2^4$&  no\cr
 2 & 270& $\{t_{1,y,z,u}\colon y,z,u\in \K^\times\}\cong \K^{\times3}$& yes\cr
 3 & 4& $\{t_{x,y,x^{-3}y^{-2},x}\colon x,y\in \K^\times\}\cong \K^{\times2}$&  no\cr
 3 & 59& $\{t_{x,y,xy,y }\colon x^3=y^3=1 \}\cong\Z_3^2$&  no\cr
 3 & 96& $\{t_{ 1,x,y,\frac1{xy}}\colon x,y\in \K^\times\}\cong \K^{\times2}$& yes\cr
 4 & 2& $\{t_{x,\frac{y}{x^2},x,x}\colon y^2=1\}\cong \K^{\times }\times\Z_2$&  no\cr
 4 & 7& $\{t_{x,\frac{y}x,1,1}\colon y^2=1\}\cong \K^{\times}\times\Z_2$&  no\cr
 4 & 30& $\{t_{1,y,z,\frac1{y^2z}}\colon       y,z\in \K^\times\}\cong \K^{\times2}$&  no\cr
 4 & 34& $\{t_{x,y,x,x}\colon x^2=y^2=1\}\cong \Z_2^2$&  no\cr
 4 & 46& $\{t_{y^2,y,y^2u,u}\colon u^2=1=y^4\}\cong \Z_2\times \Z_4$&  no\cr
 6 & 10& $\{t_{1,x,\frac1{x^2},1 }\colon x\in \K^\times\}\cong \K^{\times}$&  no\cr
 6 & 11& $\{t_{x^{-\frac23},x,x^{-\frac23},1 }\colon x\in \K^\times\}\cong \K^{\times}$&  no\cr
 6 & 20& $\{t_{x,y,x,x }\colon x^2=y^2=1\}\cong \Z_2^2$&  no\cr
 6 & 55& $\{t_{x,y,1,y }\colon x^2=y^2=1\}\cong \Z_2^2$&  no\cr
 6 & 56& $\{t_{ 1,x,\frac1{x^2},x}\colon x\in \K^\times\}\cong \K^{\times}$&  no\cr
 6 & 78& $\{t_{1,1,x,\frac1x }\colon x\in \K^\times\}\cong \K^{\times}$&  no\cr
 6 & 318& $\{t_{1,1,1,1 } \}$&  no\cr
 8 & 8& $\{t_{1,y,1,1}\colon y^2=1\}\cong \Z_2$&  no\cr
 12 & 58& $\{t_{1,1,1,1 } \}$&  no\cr
\hline
\end{tabular}
\end{center}

\begin{center}
\textsc{Table 2}
\end{center}

 Some comments on the torality of the above quasitori follow.
 The torality of $Q(96,\id)$ and $Q(270,\id)$ can be obtained by applying
 Lemma~\ref{toroporalgotoral} below, since $\si_{96},\si_{270}\in\W$.
   The remaining cases are nontoral: $Q(318,\id)$  because
   $\si_{318}\notin\W$,
 and for any other index $j$ because the identity component $L_e$ in the grading induced by $Q(j,\id)$
 verifies that $\dim L_e<4$ (as in \cite[\S~2.4]{g2} the torality of
 a grading can be characterized by $\text{rank}\,L_e=\text{rank}\, L$, taking into account
   that $L_e$ is a reductive Lie algebra.)

The problem here is the proliferation of nontoral quasitori.
Observe that $Q(j,\id)$ is nontoral for $22$ of the $25$ chosen
indices. This   is not surprising, because the group $\aut\d4$ is
smaller than $\aut\f4$ but of the same rank (in fact, $\V$ is
isomorphic to the Weyl group of $\f4$), thus there are not so many
elements to conjugate. That is why our aim will be finding only
the maximal quasitori.

Notice that the Lie algebra fixed by $\tilde\sigma_3$ is of type
$\mathfrak b_3$, so $\tilde\sigma_3$ is one of the order $2$
automorphisms providing the symmetric pair $(\d4,\mathfrak b_3)$.
Concretely $s_6$ is in the orbit of $\sigma_3$. On the other hand,
$s_5$ is in the orbit of $\sigma_4$ so that $\tilde\sigma_4$ is an
order $3$ automorphism fixing a Lie algebra of type $\mathfrak
g_2$. Thus $\aut\d4=\Int\d4\cdot\span{1,\tilde\sigma_3}
\cdot\span{1,\tilde\sigma_4,\tilde\sigma_4^2}$. In fact,
$Q(3,\id)=Q_4\cong\Z_2\times(\K^\times)^3$ and
$Q(4,\id)=Q_{13}\cong\Z_3\times(\K^\times)^2$ are maximal
quasitori, since the automorphisms commuting with $\tilde
\sigma_3$ (respectively, $\tilde \sigma_4$) are just the
extensions of automorphisms of $\mathfrak b_3$ (respectively of
$\g2$), so that if we consider the maximal torus of $\mathfrak
b_3$ (resp. $\g2$) we obtain $Q_4$ (resp. $Q_{13}$.) This will be
explained with more detail in the next section.

In the proof of Theorem~\ref{losmads} and along the paper, we use extensively the
following technical results.

\begin{lm}\label{subconjugar}\cite[Prop.~7,\,p.~27]{f4}
Let $F=\{f_{0},f_1,\ldots,f_n\}\subset\aut\d4$ be a nontoral
commutative family of semisimple elements such that
$\{f_1,\ldots,f_n\}\subset\Int\d4$ is toral.  Then,   the subgroup
generated by $F$ is conjugated to some subgroup of the form
$\span{f,t_1,\ldots,t_n}$ where $t_i\in \To$ and $f \in\No$ is
 conjugated to $f_{0}$. Moreover, this can be done in such a way that $F\cap\To\subset
\span{f,t_1,\ldots,t_n}$.
\end{lm}

\begin{lm}\label{toroporalgotoral}
If $L$ is a simple Lie algebra, $T$ is a torus of $\aut L$ and $H$ is a toral subgroup of
$\aut L$ commuting with $T$, then $HT$ is toral.
\end{lm}
\textbf{Proof.} Let $Z$ be the centralizer of $H$ in $\aut L$. As
$H$ is toral, there is $T'$ a maximal torus of $\aut L$ such that
$H\subset T'$. Hence $T'\subset Z$ and it is also a maximal torus
of $Z$. But $T\subset Z$ so that  there is $p\in Z$ such that
$pTp^{-1}\subset T'$. Consequently $p(HT)p^{-1}=HpTp^{-1}\subset
HT'\subset (T')^2\subset T'$ and $HT$ is contained in the torus
$p^{-1}T'p$.
 $\square$\smallskip

 \begin{lm}\label{moviendo}
 Fix $t\in\To$ and $j\in\{1,\dots,1152\}$.
 \begin{itemize}
 \item If $\tor{j}$ is finite, there is $s\in\To$ such that $\Ad s(\tilde\si_jt)=\tilde\si_j$.
 \item If $\tor{j}$ is not finite, there are $s\in\To$ and $t'\in\tor{j}$ such that $\Ad s(\tilde\si_jt)=\tilde\si_jt'$.
 \end{itemize}
 Therefore $Q(j,t)$ is conjugated to $Q(j,\id)$ in all the cases.
 In particular  $Q_{12}\cong Q(20,t)$,
$Q_{13}\cong Q(4,t)$ and $Q_{14}\cong Q(59,t)$ for any $t\in\To$.
 \end{lm}

 \textbf{Proof.} According to the proof and notations of \cite[Prop.~6,\,p.~26]{f4}, it is
 enough to check that $\tor{j}\cap\sor{j}$ is finite for   the indices $j$ corresponding to the representatives
 of the orbits. This is a straightforward computation.  For instance,
in our remarked cases $\tor{59}\cong\Z_3^2$,
$\tor{20}\cong\Z_2^2$, and, although $\tor{4} $ is not finite,
$\tor{4}\cap\sor{4}=\{t_{x,1,1,x}\colon x^3=1\}\cong\Z_3$ is so.
$\square$

Notice for further use that the conjugation automorphism
has been taken $\Ad s$ for certain $s\in\To$, so that it
does not \lq\lq move" the normalizer $\No$ nor any element in $\To$.\smallskip

%\begin{lm}\label{Kasper}\cite[Lemma~1.1.3, p.\,5]{Kasper}
%Every subquasitorus $Q$ of $\Int\d4$ such that $\X(Q)$ has two
%generators is toral.
%\end{lm}
%\textbf{Proof.} Following \cite[Lemma\,2]{f4}, it is enough to note that that $\Int\d4$
%(the connected component of $\aut\d4$)
%is a connected reductive group whose derivated subgroup is simply
%connected. $\square$ \margen{falso}

\begin{re}\label{casiKasper} \rm One of the most useful tools   in the computational approach to the gradings
on $\f4$   was  \cite[Lemma~2]{f4}, according to which every quasitorus of $\aut\f4$ such that
$\X(Q)$ has two generators is toral. Of course we cannot apply this result to $\aut\d4$,
 because one single outer automorphism
  generates a nontoral quasitorus. But even if we consider a subquasitorus $Q$ of the connected component $\Int\d4$, it could happen that $Q$ were nontoral with two generators. Indeed, take
$Q=\langle\tilde\si_1,t_{-1,1,1,1}\rangle\cong\Z_2^2$, which is
nontoral but $\si_1\in\W$ (equivalently, $\tilde\si_1$ is an inner
automorphism.) The nontorality of $Q$ is consequence of the
nontorality of $Q(1,\id)\cong Q\times(\K^\times)^2$, by applying
Lemma~\ref{toroporalgotoral}. Moreover, it is possible to prove
that if $Q$ is a nontoral subquasitorus of $\Int\d4$ such that
$\X(Q)$ has two generators, $Q$ can be subconjugated inside
$Q(1,\id)$.  The source of these differences between the behavior
of $\aut\f4$ and $\Int\d4$ is that $\Int\d4$  is not simply
connected. Anyway, what we will use is a much weaker result: Note
that if   $Q$ is the quasitorus generated by $
\{\tilde\si_it_1,t_2\}$  with $\si_i\in\W$ of order 3 and
$t_j\in\To$, then $Q$ is obviously toral by
Lemma~\ref{toroporalgotoral} ($\si_i$ would be conjugated to
$\si_{96}$), and also if $\si_i\in\W$ is of any other order but
$t_2$ is contained in a subtorus of $\tor{i}$.
\end{re}

  \textbf{Proof of Theorem~\ref{losmads}.} First of all notice that these
quasitori $Q_i$ are not   conjugated, since the grading groups
$\X(Q_i)$ are not isomorphic. Besides, we have
remarked that $Q_1,\dots,Q_{11}$ are conjugated to the MAD-groups
obtained from \cite{LGII}, which cover all the MAD's contained in
$\Int\d4\cdot\span{1,\tilde\sigma_3}$. Thus, what we are going to
prove is that $Q_{12},Q_{13}$ and $Q_{14}$ are maximal quasitori,
and that if $Q$ is a maximal quasitorus not conjugated to any
subgroup of $\Int\d4\cdot\span{1,\tilde\sigma_3}$, then it is
conjugated to one of these three quasitori.

We now check  that  $Q_{14}$ is a maximal quasitorus.  Let us
prove that if $f\in\aut\d4$ commutes with $Q_{14}$, then $f$
  belongs to $Q_{14}$. With that purpose consider
$Z=\Co_{\aut\d4}(\tor{59})$. The group $\span{Q_{14}\cup\{f\}}$ is
an abelian subgroup of $Z$, as well as its closure   in the Zarisky
topology. But this is again a quasitorus, whence it is contained
in the normalizer of some maximal torus $T$ of $Z$. In particular
$\span{Q_{14}\cup\{f\}}\subset N_Z(T)$. By construction also
$\To\subset Z$ so that there is some $p\in Z$ such that
$pTp^{-1}=\To$. Consequently $p\span{Q_{14}\cup\{f\}}p^{-1}\subset
 N_{\aut\d4}(\To)=\No$ and $ pfp^{-1}, p\tilde\si_{59}p^{-1}\in
 \No\cap\Co_{\aut\d4}(\tor{59}) $ with $ptp^{-1}=t$ for any
$t\in\tor{59}$. Take $j_1,j_2\in 1,\dots,1152$ and $t_1,t_2\in\To$
such that
$pfp^{-1}=\tilde\sigma_{j_1}t_1,p\tilde\sigma_{59}p^{-1}=\tilde\sigma_{j_2}t_2$.
A straightforward computation in the computer tell us that $\{i\in
1,\ldots,1152\colon \sigma_i\cdot
t_{1,\omega,\omega,\omega}=t_{1,\omega,\omega,\omega},\sigma_i\cdot
t_{\omega,1,\omega,1}=t_{\omega,1,\omega,1} \}=\{59,835,894\}$. In
particular   $j_i$ must be one of those indices, and
$\sigma_{j_i}=\sigma_{59}^{n_i}$ for $n_i\in\{0,1,2\}$. Moreover,
$n_2\ne0$, otherwise  the grading produced by $pQ_{14}p^{-1}$
would be toral. We can take $t_2=\id$ by replacing $p$ by $sp$,
being $s$ the element in $\To$ such that
$s\tilde\sigma_{59}^{n_2}t_2s^{-1}=\tilde\sigma_{59}^{n_2}$ as in
Lemma~\ref{moviendo}. As $\tilde\sigma_{59}^{n_1}t_1$ and
$\tilde\sigma_{59}^{n_2}$ commute, as well as
$\tilde\sigma_{59}^{n_1}$ and $\tilde\sigma_{59}^{n_2}$, then
$t_1$   commutes with $\tilde\sigma_{59}$ so that   $t_1\in
\tor{59}$. So $pfp^{-1}=\tilde\sigma_{59}^{n_1}t_1\in Q(59,\id)
=pQ_{14}p^{-1}$ and $f\in Q_{14} $.

Now let  $f\in C_{\aut\d4}(Q_{13})$  and as before find an
automorphism $p\in\aut\d4$ such that $ptp^{-1}=t$ for all
$t\in\tor{4}$ and $pfp^{-1}=\tilde\sigma_{j_1}t_1,\,
p\tilde\si_{4}p^{-1}=\tilde\sigma_{j_2}t_2\in C_{\No}(
\tor{4})$
for $t_1,t_2\in\To$. Thus
$\sigma_{j_i}\in\{\sigma_k\in\V\colon\tor{4}\subset\tor{k}\}=
\{\sigma_4^n\colon n=0,1,2\}\cdot\{  1,\sigma_3\}$ (isomorphic to
the group of permutations $S_3$.) Besides $ p\tilde\sigma_4p^{-1}$ has
order 3, so there is $l\in\{1,2\}$  such that $\si_{j_2}=\si_4^l$.
Now  we can assume that $t_2\in\tor{4}$, by replacing $p$ by $sp$,
being $s$ the element as in Lemma~\ref{moviendo} verifying
$s\tilde\sigma_{4}^{l}t_2s^{-1}\in Q(4,\id)$.
 Thus $pQ_{13}p^{-1}=Q(4,\id)$. Note
that $\sigma_4$ commutes neither with $\sigma_3$ nor
$\sigma_3\sigma_4$ nor $\sigma_3\sigma_4^2$, since
$\sigma_3\sigma_4=\sigma_4^2\sigma_3$ with $\sigma_3$ of order
$2$. But $pfp^{-1}$ commutes with $p\tilde\sigma_4p^{-1}$, so
their projections by $\pi$ commute and
$\sigma_{j_1}\in\{\sigma_4^n\colon n=0,1,2\}$. In particular
$\tilde\si_{j_1}$ commutes with  $\tilde\sigma_{4}^{l}t_2$ and
hence $t_1$ does so. Hence $t_1 \in \tor{4}$, $pfp^{-1}
=\tilde\sigma_{j_1} t_1\in Q(4,\id)=pQ_{13}p^{-1}$ and $f\in
Q_{13}$.

Finally, for  $f\in C_{\aut\d4}(Q_{12})$  we again find an
automorphism $p\in\aut\d4$ such that $ptp^{-1}=t$ for all
$t\in\tor{20}$ and $pfp^{-1} ,\, p\tilde\si_{20}p^{-1}\in C_{\No}(
\tor{20})$. Denote
${J}=\{\sigma_k\in\V\colon\tor{20}\subset\tor{k}\}$. As
$p\tilde\si_{20}p^{-1}=\tilde\sigma_jt$ for some $\sigma_j\in J$,
in particular $\sigma_j^2\in{J}$ has order exactly $3$ (not $1$
because the grading induced by $
\span{\tilde\sigma_{20}^2,\tor{20}}$ is also nontoral.)
 But there are
$32$ order 3 elements in ${J}$, all of them in the orbit of
$\sigma_4$. Hence $(\tilde\sigma_jt)^2$ is conjugated to
$\tilde\sigma_4t_1$ without moving $\No$, but also without moving
$\tor{20}$, since the corresponding $32$ conjugating elements in
$\mathcal{V}$ can be taken inside ${J}$. Selecting now the order
$2$ elements $\sigma_i\in J$
   commuting with $\sigma_4$ (looking for $\sigma_j^3$), there is only one,
$\sigma_{259}=-\id$ (here $\id$ denotes the $4\times4$ identity matrix.) Replacing $p$, we have found that
$pQ_{12}p^{-1}=\span{\tilde\sigma_4t_1,\tilde\sigma_{259}t_2,\tor{20}}=:Q$
with $pfp^{-1}\in\No$ commuting with $Q$. Notice now that
$\tilde\sigma_{40}Q\tilde\sigma_{40}^{-1}=Q(20,t_3)$, so there is
$s\in\To$ such that
$(s\tilde\sigma_{40})Q(s\tilde\sigma_{40})^{-1}=Q_{12}$. Replace
again $p$ by $s\tilde\sigma_{40}p$  to get $pfp^{-1}\in\No$
commuting with $pQ_{12}p^{-1}=Q_{12}$. As there are only $6$
elements $\sigma_i\in J$ with
$\sigma_i\sigma_{20}=\sigma_{20}\sigma_i$, they must be just
$\{\sigma_{20}^n\colon n=0,\dots 5\}$. Thus there is $n$ such that
$pfp^{-1}=\tilde\sigma_{20}^nt_4$, which commutes with
$\tilde\sigma_{20}\in pQ_{12}p^{-1}$, so that $t_4\in\tor{20}$ and
$pfp^{-1}\in Q_{12}=pQ_{12}p^{-1}$. We have already  seen that
$Q_{12}$,  $Q_{13}$ and  $Q_{14}$ are maximal abelian diagonalizable groups.

Suppose now that we have another $Q$ maximal quasitorus contained
in $\No$ such that nor $Q$ neither any of its conjugated quasitori
are contained in $\Int\d4\cdot\{1,\tilde\sigma_3\}$. Recall that
there are three orbits (under conjugation in
$\mathcal{V}$) of order 3 elements in $\mathcal V$, given by the representatives
$\sigma_4$, $\sigma_{59}$ and $\sigma_{96}$.  We are going to
prove that the abelian group $\pi(Q)$ intersects   the orbit of
either $ \sigma_{59}$ or $ \sigma_{4}$. Thus, we will be able to
assume that either $\tilde\sigma_{59}t$ or $\tilde\sigma_{4}t$
belongs to $Q$ (for some $t\in\To$), by conjugating by an
element in $\No$. Indeed, consider $\mathcal{V}_0=\mathcal{W}$,
$\mathcal{V}_1=\mathcal{W}\sigma_3$,
$\mathcal{V}_2=\mathcal{W}\sigma_3\sigma_4$, $\mathcal{V}_3
=\mathcal{W}\sigma_3\sigma_4^2$, $\mathcal{V}_4
=\mathcal{W}\sigma_4$ and $\mathcal{V}_5=\mathcal{W}\sigma_4^2$
(obviously $\V=\cup_{i=0}^6\V_i$.) By hypothesis, $\pi(Q)$ is not
contained in $\mathcal{V}_0\cup\mathcal{V}_1$. It is neither
contained in $\mathcal{V}_0\cup\mathcal{V}_i$, for $i=2,3$, since
they are conjugated to $\mathcal{V}_0\cup\mathcal{V}_1$ by means
of $\sigma_4$ and $\sigma_4^2$ respectively, since
$\sigma_3\sigma_4=\sigma_4^2\sigma_3$. Moreover, $\pi(Q)$ is not
contained in
$\mathcal{V}_0\cup\mathcal{V}_1\cup\mathcal{V}_2\cup\mathcal{V}_3$,
since $\pi(Q)$ is abelian, and none of the elements in $\mathcal{V}_i$
commutes with none of the elements in $\mathcal{V}_j$ nor
$\mathcal{V}_k$, for $\{i,j,k\}$ distinct indices in $\{1,2,3\}$.
Hence $\pi(Q)\cap(\mathcal{V}_4\cup\mathcal{V}_5)\ne\emptyset$,
but every element in $\mathcal{V}_4\cup\mathcal{V}_5$ has order
multiple of $3$ ($3,6,12$) and one of its powers is an order 3
element in the orbit   of
$\sigma_{59}$ or $\sigma_4$ ($\sigma_{96}\in\mathcal{W}$), which
also belongs to $\pi(Q)$.

First, suppose that $\si_4\in\pi(Q)$. Note that $\pi(Q)=
(\pi(Q)\cap(\V_0\cup\V_1))\cdot\langle\si_4\rangle$, but, since it
is abelian, $\pi(Q)= \pi(Q)\cap\W\cdot\langle\si_4\rangle$ (take
into account that $\si_4$ does not commute with any element in
$\V_{1,2,3}$.) Suppose that there is an order 3 element in
$\pi(Q)\cap\W$. As $\{\si\in\W\colon \si\si_4=
\si_4\si,\si^3=\id\}=\langle\si_{952}\rangle\cong\Z_3$, we have
that $\si_{952}\in\pi(Q)$, but $\{\si\in\W\colon \si\si_4=
\si_4\si,\si\si_{952}=\si_{952}\si\}=\langle\si_{952},\si_{259}=-\id\rangle\cong\Z_6$,
so that either $\pi(Q)=\langle\si_{952},\si_4\rangle$ or
$\pi(Q)=\langle\si_{952},-\si_4\rangle$. In the first case there
are $t_1,t_2\in\To$ with $Q=\langle
\tilde\sigma_{4}t_1,\tilde\sigma_{952}t_2, t_3\rangle$, for
$\langle t_3\rangle=\tor{4}\cap\tor{952}=\{t_{x,1,1,x}\colon
x^3=1\}\cong\Z_3$. As $\langle \tilde\sigma_{952 }t_2,t_3\rangle$
is necessarily toral by Remark~\ref{casiKasper}, apply
Lemma~\ref{subconjugar} to get that $Q$ is conjugated to $\langle
t_3,t_4,\tilde\sigma_{i}t_5\rangle$ with $\tilde\sigma_{i}t_5$
conjugated to $\tilde\sigma_{4}t_1$, and, in particular,
$\si_i\notin\W$ ($\si_i$ could be in other orbit, and could have
order 3, 6 and 12, but still belongs to $\V_4\cup\V_5$.) By
maximality, $Q$ is conjugated to $Q(i,t_5)\cong Q(i,\id)$ and
$\tor{i}=\langle t_3,t_4\rangle\cong\Z_3\times\Z_r$. The only
possibility, according to Table~2, is that $\si_i$ is  in the
orbit of $\si_{59}$ and $Q$ is conjugated to $Q(59,\id)=Q_{14}$.
In the second case, there are $t_1,t_2 \in\To$ with $Q=\langle
\tilde\sigma_{1149}t_1,\tilde\sigma_{952}t_2\rangle$
($-\si_4=\si_{1149}$), since
$\tor{1149}\cap\tor{952}=\{t_{1,1,1,1}\}$. Obviously  the grading
induced by $\langle \tilde\sigma_{952}t_2\rangle$ is   toral and,
by  Lemma~\ref{subconjugar} there are $t_4,t_5\in\To$ such that
$Q$ is conjugated to $\langle t_4,\tilde\sigma_{i}t_5\rangle$ for
$t_4$ of order $3l$. By  maximality of $Q$, this set coincides
with $Q(i,t_5)$ and hence $\tor{i}\cong\Z_{3l}$, what it is a
contradiction. Now we suppose that $\pi(Q)\cap\W$ does not contain
order 3 elements. Hence $\pi(Q)\cap\W\subset\{\pm\id,\pm\si_{111},
\pm\si_{211},\pm\si_{249}\}$ (the only elements in $\W$ commuting
with $\si_4$ of order coprime with 3.) Consequently the
possibilities for $\pi(Q)$ are: $ \langle\si_4\rangle$, $
\langle-\si_4\rangle$, $\ \langle\si_4,\si_j\rangle$ and $
\langle-\si_4,\si_j\rangle$, for some $j\in\{111,211,249\}$. If
$\pi(Q)=\langle\si_4\rangle$, then $Q=Q(4,t)\cong
Q(4,\id)=Q_{13}$. As $-\si_4$ is in the orbit of $\si_{20}$, if
$\pi(Q)=\langle-\si_4\rangle$, then $Q\cong Q(20,t)\cong
Q(20,\id)=Q_{12}$. If $\pi(Q)=\langle\si_4,\si_j\rangle$, there
are $t_1,t_2\in\To$ such that
$Q=\langle\tilde\si_4t_1,\tilde\si_jt_2\rangle\cdot\tor{4}\cap\tor{j}$.
But $\tor{4}\cap\tor{j}\cong\K^\times$ for any of the three
indices, so the grading induced by $
\{\tilde\si_jt_2\}\cup(\tor{4}\cap\tor{j})$ is toral by
Remark~\ref{casiKasper}, and by Lemma~\ref{subconjugar} there are
$t_3,t_4\in\To$ and an index $i$  such that
$Q=\langle\tilde\si_it_3,t_4\rangle\cdot\tor{4}\cap\tor{j}=Q(i,t_3)$.
In particular $\tor{i}=\langle
t_4\rangle\cdot\tor{4}\cap\tor{j}\cong\Z_{2s}\times\K^\times$, a
contradiction (the only possibilities  with
$\tor{i}\cong\Z_{2s}\times\K^\times$ would be $i=2,7$, up to
conjugation, but $\si_7\in\V_0$ and $\si_2\in\V_2$.) Finally, if
$\pi(Q)=\langle-\si_4,\si_j\rangle$, there are $t_1,t_2\in\To$
such that $Q=\langle
\tilde\si_{1149}t_1,\tilde\si_jt_2\rangle\cdot\tor{1149}\cap\tor{j}$.
 Although $\tor{1149}\cap\tor{j}\cong\Z_2$, we can apply  Remark~\ref{casiKasper} to get
 that the grading induced by $ \{ \tilde\si_jt_2\}\cup(\tor{1149}\cap\tor{j})$ is toral,
 because $\tor{1149}\cap\tor{j}\subset\tor{4}\cap\tor{j}\cong\K^\times$. Thus,
by Lemma~\ref{subconjugar}, the quasitorus $Q$ is conjugated to
$Q(i,t)$ for certain index $i$ verifying
$\tor{i}\cong\Z_{2l}\times\Z_2$. The only indices of
representatives in these conditions are 34, 46, 20 and 55, but
$\si_{55},\si_{34}\in\V_0$ and $\si_{46}\in\V_1$ so that $\si_i$
is in the orbit of $\si_{20}$, and $Q\cong Q(20,\id)=Q_{12}$.

 Second,
suppose that $ \sigma_{59}\in \pi(Q)$. This time  we can suppose
that $\tilde\sigma_{59}\in Q$, by Lemma~\ref{moviendo}. As before,
$\pi(Q)= (\pi(Q)\cap\W)\cdot\langle\si_{59}\rangle$ and by
maximality $Q\cap\To=
\cap\{\tor{i}\cap\tor{59}\colon\si_i\in\pi(Q)\cap\W\}$. In
particular $Q\cap\To\subset\tor{59}\cong\Z_3^2$, so that
$Q\cap\To\cong\Z_3^{0,1,2}$. If $Q\cap\To=\tor{59}$, then
$Q_{14}=Q(59,\id)\subset Q$ and, by maximality of $Q_{14}$ it
follows $Q=Q_{14}$. If $Q\cap\To=\span{t_1}\cong\Z_3$, then
$\{\id\}\ne\pi(Q)\cap\W\subset\{\sigma_i\in\W\colon
\sigma_i\sigma_{59}=\sigma_{59}\sigma_i\}$, which is a subgroup of
$\V$   with only one order $2$ element ($-\id=\sigma_{259}$.) If
there were some element in $\pi(Q)$ of order not divisor of $3$,
then $\sigma_{259}\in\pi(Q)$ and $\Z_3\cong
Q\cap\To\subset\tor{259}\cong\Z_2^4$, what is absurd. Hence
$\pi(Q)\subset S:= \{\sigma\in\W\colon
\sigma^3=\id,\sigma\sigma_{59}=\sigma_{59}\sigma\}$, which is a
non abelian set with $9$ elements. There must be $\sigma_i\in S$
such that $\pi(Q)\cap\W=\span{\sigma_i}$ (otherwise another
$\si_j\in S \setminus\langle\si_i \rangle$ would satisfy
$\langle\si_i,\si_j\rangle\subset\pi(Q)\cap\W$ but then $
\pi(Q)\cap\W$ would have at least 9 elements  belonging to $S$,
the whole $S$, but $S$ is not abelian.) Thus, there is some
$t_2\in\To$ such that
   $ \span{\tilde\sigma_{59},\tilde\sigma_{i}t_2,t_1}=Q$.   As
$\{\tilde\sigma_{i}t_2,t_1\}$ induces a toral grading   by
Remark~\ref{casiKasper} ($\si_i$ has order 3), we can conjugate
$Q$ to $\span {t_1,t_3,f=\tilde\sigma_{k}t_4}$ with $f$ a
conjugate of $\tilde\sigma_{59} $ (hence of order just 3.) We can
take $\sigma_{k}=\sigma_{59}$ or $\sigma_{4}$ by conjugating  now
inside $\No$. In the first case $Q=Q(59,t_4)\cong Q_{14}$, and in
the second case $Q\subsetneqq Q(4,t_4)$, a contradiction with the
maximality of $Q$. Finally suppose that $Q\cap\To=\{\id\}$. If
$\pi(Q)\cap\W\subset S$, $\pi(Q)\cap\W$ would have only one order
3 generator    and $Q$ would be contained strictly in one of the
quasitori in the paragraph above. Thus
$-\id=\sigma_{259}\in\pi(Q)$. Notice also that
 $\span{\tilde\sigma_{59},
  \tilde\sigma_{259}t }$ is not a MAD (the second automorphism is inner, so the set is
conjugated by Lemma~\ref{subconjugar} to $\span{\tilde\sigma_{k}t_1,t_2}$, where $\tilde\sigma_{k}t_1$
has order 3, and so $\si_k$ is in the orbit of $\si_{59}$ or $\si_{4}$, hence, it is strictly contained
in $Q(4,t_1)$ or $Q(59,t_1)$), so that there is
   $\sigma_i\in\pi(Q)\cap\W\setminus\span{ \sigma_{259}}$.
Let us check that $\si_i$ cannot have order 4. In such case,
 $\sigma_i\si_{59}$ would have order 12 and it would be  conjugated to $\sigma_{58}$, so there would be
$t\in\To$ such that $A=\span{(\tilde\sigma_{58}t)^n\colon
n=1\dots12}\subset Q\subset\No$ (up to conjugation), but the centralizer
$C_{\No}(A)=A$ (since
$\tor{58}=\{\id\}$ and   $\{\sigma\in\V\colon
\sigma\sigma_{58}=\sigma_{58}\sigma\}$ has 12 elements, obviously
the powers of $\sigma_{58}$), so that $Q=A$. This is an absurd
since $A$ is not maximal in $\aut\d4$ (it is only maximal in $\No$.)
Because $\sigma_{58}^3\in\W$, so
$A\cong\span{(\tilde\sigma_{58}t)^4,(\tilde\sigma_{58}t)^3}$ is
conjugated to $\span{\tilde\sigma_{k}t_1,t_2}$, strictly contained
in $Q(4,t_1)$ or $Q(59,t_1)$, as above.
Therefore $\sigma_i$  has order either 6 or 3.
It can be taken of order $3$ (if $\sigma_i^3=-\id$
 then $-\sigma_i\in\pi(Q)$ has order 3.)
  Then $\pi(Q)=\span{\sigma_{59},\sigma_{259}\sigma_{i}}$ (any   element in $S$ different than $(\si_i)^{0,1,2}$ does not
  commute with $\si_i$),
  with $\sigma_i\sigma_{259}\in\W$ of order 6, so applying Lemma~\ref{subconjugar}
  we obtain that $Q$ is conjugated to $\span{\tilde\sigma_{k}t_1,t_2}$ with $k=4$
  ($\sigma_k\notin\W$, $\sigma_k^3=\id$ and $\Z_6\subset\tor{k}$), strictly contained
in $Q(4,t_1)$, a contradiction.
$\square$

\section{Algebraic approach}

 We are going to revisit the MAD-groups of $\aut\d4$ which intersect
 some connected component with   order 3 outer automorphisms, looking at them under the light of triality and related
 stuff. While Section 3 is of a computational nature, this is more conceptual and qualitative.

According to \cite{Kac}, there are two order 3 outer
automorphisms of $\d4$, up to conjugation, related to the affine
diagrams \vskip .5cm \hskip 2cm
 \xy <1cm, 0cm>:
 \POS(0,0)\xycircle(.1,.1){-}
 \POS(0,.1)\ar @{-}+(1,0)*+
 \xycircle(.1,.1){-}
 \POS(1.22,.1)\ar @3{-}+(1,0)*+
 \xycircle(.1,.1){-}
 \POS(1.6,.1)\ar @{-}+(.15,.15)
 \POS(1.6,.1)\ar @{-}+(.15,-.15)
 \POS(-.4,-.1)\ar @{-}+(.4,.4)
 \POS(-.4,.3)\ar @{-}+(.4,-.4)
\endxy
\vskip -.45cm \hskip 7cm
 \xy <1cm, 0cm>:
 \POS(0,0)\xycircle(.1,.1){-}
 \POS(0,.1)\ar @{-}+(1,0)*+
 \xycircle(.1,.1){-}
 \POS(1.22,.1)\ar @3{-}+(1,0)*+
 \xycircle(.1,.1){-}
 \POS(1.6,.1)\ar @{-}+(.15,.15)
 \POS(1.6,.1)\ar @{-}+(.15,-.15)
 \POS(2.05,-.1)\ar @{-}+(.4,.4)
 \POS(2.05,.3)\ar @{-}+(.4,-.4)
\endxy
\vskip .7cm

\noindent obtained from \cite[TABLE Aff\,3, p.\,55]{Kac}, which
fix subalgebras of type $\g2$ and $\mathfrak a_2$, respectively.
The first one is known by its relationship with the triality
phenomenon, but the second one is also related to it and both
cases can be treated in a uniform way.

\subsection{Automorphisms of $\d4$}

 Denote   by $(C,q)$ a
Cayley algebra with standard quadratic map $q$. As usual its
canonical involution will be denoted by $x\mapsto\bar x$. Under
our assumptions on the ground field, such a $C$ is unique up to
isomorphism. Now an algebra of type $\d4$ is
$\mathfrak{o}(C,q)=\{g\in\mathfrak{gl}(C)\colon
b(g(x),y)+b(x,g(y))=0\quad\forall x,y\in C\}$, where $b$ denotes
the polar form of $q$. In this subsection we would like to get a
description of $\mathcal G:=\aut(\d4)$ which could be useful in
our further study of quasitori and gradings on $\d4$. Denote by
$\go(C,q)$ the algebraic group of all bijective linear maps
$f\colon C\to C$ such that there is $\lambda\in \K^\times$ such
that $q(f(x))=\lambda q(x)$ for any $x\in C$. The scalar $\lambda$
is called the {\em multiplier} of $f$ and it is usually denoted by
$\mu(f)$ (see \cite[\S 12]{Boi} for more information.)
 It is well known (for instance \cite[Exercise 15, p.~287]{Jac})
 that the identity component $\mathcal G_0$ of the
algebraic group $\mathcal G$ is the group of all automorphisms of
the form $g\mapsto fgf^{-1}$ where $f\in\go(C,q)^+$ (notation as
in \cite[p.~154]{Boi}.) But two elements $f_1,f_2\in\go(C,q)^+$
induce the same automorphism $g\mapsto f_igf_i^{-1}$ ($i=1,2$) if
and only if $f_1=\lambda f_2$ for some $\lambda\in \K^\times$.
Thus we can identify $\mathcal G_0$ with
$\pgo(C,q)^+=\go(C,q)^+/\K^\times I$.

It is also well known that $\mathcal G$ is an extension $\mathcal
G=\mathcal G_0\cdot S_3$ where $S_3$ denotes the symmetric group
of order 3. The aim of this subsection is to clarify the
nature of this extension, which will enable us to make some
explicit computations (for instance on centralizers.) To describe
the extension $\mathcal G_0 \cdot S_3$ we need to explicit the
action of the generators of $S_3$ on $\mathcal G_0$. If we take an
element $f\in\go(C,q)^+$, we shall denote by $\bar f$ the new
element in $\go(C,q)^+$ such that $\bar f(x):=\ov{f(\bar x)}$ for
any $x\in C$.  The map $f\mapsto\bar f$ induces an order 2
automorphism $\sigma:\pgo(C,q)^+\to\pgo(C,q)^+$. This will be
identified with an order 2 permutation in $S_3$.  Moreover, to
describe the action on $\mathcal G_0$ of the cyclic order 3
permutation in $S_3$ we need to take into account the two (unique
up to isomorphism) possible $8$-dimensional symmetric compositions
which can be constructed from $(C,q)$ (\cite[\S~34]{Boi}.) The
first one is the para-Hurwitz algebra $\bar C$ with multiplication
$x\cdot y=\bar x\ \bar y$, while the second one is the Okubo
algebra $\oku$ whose multiplication is $x\cdot y=\varphi(\bar
x)\varphi^{-1}(\bar y)$, where $\varphi\in\aut(C)$ is an order
three automorphism fixing a four-dimensional algebra (such an
automorphism is unique up to conjugation.) It is easy to check
that $\aut(\bar C)=\aut(C)=G_2$ while $\aut(\oku)$ is the identity
connected component of $A_2=\aut(\mathfrak{sl}(3,\K))$. Though
this is a well known fact (\cite{Alb3}),  one can check it by
taking into account
 that the Okubo algebra is isomorphic to
$(\mathbb P,*)=(\text{Mat}_{3}(\K)_0,*)$ the pseudo-octonion
algebra   constructed on the vector space of zero trace $3\times
3$ matrices with entries in $\K$ (\cite{okubo, Alb3}.) Its product
is given by $x\ast y:=\mu xy+(1-\mu)yx- \frac{1}{3}\hbox{tr}(xy)$
for   $\mu=\frac{1-\omega}3$, with $\omega$ primitive cube root of
unity and $\hbox{tr}(\cdot)$ denotes the matrix trace. If
$p\in\GL(3,\K)$, then $x\mapsto pxp^{-1}$ is an automorphism of
$\mathbb P$. That is,
$\Int(\mathfrak{sl}(3,\K))\subset\aut(\mathbb P)$. On the other
hand,  the product $[x,y]:=xy-yx$ is $[x,y]=(2\omega-1)^{-1}(x\ast
y-y\ast x)$. Thus any element in $\aut(\mathbb P)$ is an
automorphism of the Lie algebra $\mathfrak{sl}(3,\K)$. Therefore
it is an automorphism of $\text{Mat}_3(\K)$ or the opposite of an
antiautomorphism of the same algebra. But this last possibility
does not provide an automorphism of $\mathbb P$ so $\aut(\mathbb
P)\cong\Int( \mathfrak{sl}(3,\K))=\PGL(3,\K)(=(A_2)_0)$.

If we denote by $S$ to any of the symmetric compositions
introduced in the previous paragraph, and by $*$ its product, then
it is well known (\cite[Proposition~35.4]{Boi}) that for any
$t\in\go(C,q)^+$ with multiplier $\mu(t)$ there are elements
$t^+,t^-\in\go(C,q)^+$ such that $(t,t^-,t^+)$ is an {\sl
admissible triple}, that is, $\mu(t)^{-1}t(x*y)=t^-(x)*t^+(y)$ for
any $x,y\in C$.  Moreover $t^+$ and $t^-$ are unique up to scalar
multiplication by some nonzero $\lambda$ and $\lambda^{-1}$
respectively. It is also a standard fact that $(t,t^-,t^+)$ is
admissible if and only if $(t^-,t^+,t)$ is admissible. Now the map
$t\mapsto t^-$ induces an order 3 automorphism
$\theta_S:\pgo(C,q)^+\to\pgo(C,q)^+$. We shall write
$\theta:=\theta_{\bar C}$ and $\theta':=\theta_{\oku}$. Denoting
by $[t]$ the equivalence class in $\pgo(C,q)^+$ of
$t\in\go(C,q)^+$, we have $\theta_S([t])=[t^-]$ so that
$\theta_S^3=1$ and it can be easily proved that
$\theta_S\sigma=\sigma\theta_S^2$, which together with
$\sigma^2=1$ provide a group monomorphism
$S_3\to\aut(\pgo(C,q)^+)$. Thus we get a description of
$\aut(\d4)$ as a semidirect product
$$
\mathcal{G}\cong\pgo(C,q)^+\rtimes S_3
$$
where the product of elements is given by
\begin{equation}\label{actionS3}
\sigma\cdot [t]= [\bar t\,]\sigma, \quad \theta_S\cdot [t]=[t^-]\theta_S,
\quad \theta_S^2\cdot [t]=[t^+]\theta_S^2
\end{equation}
while products of
$\theta_S$'s and $\sigma$'s are governed by the corresponding
relations in $S_3$.

Observe now that we can consider $\aut S\subset  \mathcal G_0$ independently of the
symmetric composition algebra considered $(S,*)$. In the para-Hurwitz case,
 $\aut \bar C=\aut C\subset \text{O}(C,q)^+
\subset\go(C,q)^+$, and, by composing with the canonical
projection $\go(C,q)^+\to \pgo(C,q)^+$, we still have a
monomorphism $\aut(C)\hookrightarrow \pgo(C,q)^+$. In the Okubo case, straightforward computations
show that $(x*y)*x=q(x)y$, so that any $f\in\aut(S,*)$ verifies $q(f(x))=q(x)$, that is,
$\aut S\subset  \text{O}(C,q)$. Moreover, $\aut S\subset  \text{O}(C,q)^+$ because $\aut S$ is connected ($=(A_2)_0$.)
As before, the composition of $\aut S\hookrightarrow \go(C,q)^+$ with the projection $\go(C,q)^+\to \pgo(C,q)^+$
is still an embedding (the only multiple of $f\in\aut S$ which is also an automorphism is $f$ itself.)

Now we can compute certain centralizers easily.

\begin{pr}\label{centr}
For $S=\bar C$ or $\oku$, the centralizer in $\mathcal G$ of
$\theta_S$ is $C_{\mathcal G}(\theta_S)=\aut(S)\cup \aut(S)
\theta_S\cup \aut(S)\theta_S^2$.
\end{pr}
\textbf{Proof.}  Recalling that $\mathcal G_0=\pgo(C,q)^+$ we have
$$\small \mathcal G=\mathcal G_0\cup \mathcal G_0\theta_S\cup
\mathcal G_0\theta_S^2\cup \mathcal G_0\sigma \cup \mathcal G_0
\sigma\theta_S\cup \mathcal G_0\sigma\theta_S^2.$$ Taking an
element   $[f]\in \mathcal G_0=\pgo(C,q)^+$ and imposing the
condition that it centralizes $\theta_S$, we get
$\theta_S[f]=[f^-]\theta_S=[f]\theta_S$ which implies that there
is an admissible triple $(f, f,f^+)$ and consequently an admissible triple
$(f,f^+,f)$. Since $f^-$ and $f^+$ are uniquely determined by $f$
(up to nonzero scalar multiples), we have $[f]=[f^+]=[f^-]$. This
implies that a certain nonzero multiple of $f$ is an automorphism
of $S$. So $[f]\in \aut(S)$. If we take now $n=1$ or $2$ and an
element $[f]\theta_S^n\in \mathcal G_0\theta_S^n\cap C_{\mathcal
G}(\theta_S)$, by imposing the commutativity condition we get also
$[f]=[f^+]=[f^-]$ as before. The rest of the possibilities
($[f]\in \mathcal G_0\sigma\theta_S^{0,1,2}$) are not compatible
with the commutativity condition. $\square$\smallskip

\begin{co}\label{key}
Let $S$ be the para-Hurwitz or the Okubo algebra. For any maximal
quasitorus $Q$ of $\aut(S)$ the group $\hat
Q:=\span{Q\cup\{\theta_S\}}$ is a maximal quasitorus of $\mathcal
G=\aut\d4$.
\end{co}
\textbf{Proof.}  Suppose that $r\in\mathcal G$ is a semisimple
element commuting with $\hat Q$. Then $r\in C_{\mathcal
G}(\theta_S)=H\cup H\theta_S\cup H{\theta_S}^2$, for $H:=\aut(S)$.
If $r\in H$ then $r\in C_H(Q)$ and by the maximality of $Q$ we
have $r\in Q\subset\hat Q$. In case that $r\in H\theta_S$ then
$r\theta_S^2\in H$ so that $r\theta_S^2\in C_H(Q)$ and therefore
$r\theta_S^2\in Q$ implying $r\in\hat Q$ (the other possibility
$r\in H\theta_S^2$ is similar.) $\square$\smallskip

%\margen{Se podr\'{\i}a suprimir este remark}
%\begin{re}{\rm
%It is not difficult to relate the two order three automorphisms $\theta$ and
%$\theta'$. Indeed, a triple $(t,t^-,t^+)$ is admissible relative to the
%algebra $\bar C$
%if and only if $(t,\varphi^{-1}t^-\varphi,
%\varphi t^+\varphi^{-1})$ is admissible relative to $\oku$. As a
%consequence we have: $\theta'([t])=\varphi^{-1}\theta[t]\varphi$ or
%$\theta'=\hbox{In}(\varphi^{-1})\theta$. To be more precise, we could
%consider $S_3=\langle\theta,\sigma:\theta^3=\sigma^2=1, \sigma\theta=
%\theta^2\sigma\rangle$ and the semidirect product $\pgo(C,q)^+\rtimes S_3$
%or else $S'_3=\langle\theta',\sigma:\theta'^3=\sigma^2=1, \sigma\theta'=
%\theta'^2\sigma\rangle$ and then $\pgo(C,q)^+\rtimes S'_3$. Then
%the isomorphism
%$\hbox{In}(\varphi^{-1}):\pgo(C,q)^+\to\pgo(C,q)^+$ such that $[t]\mapsto
%[\varphi^{-1}t\varphi]$ can be extended to a unique isomorphism
%$\pgo(C,q)^+\rtimes S_3\cong\pgo(C,q)^+\rtimes S'_3$ such that $\theta\mapsto
%\theta'$ and $\sigma\mapsto\sigma$. Thus we shall not distinguish both
%semidirect products in the sequel.}
%\end{re}

We know that the order 3 outer automorphisms of $\d4$ fall
into two categories: those whose fix a subalgebra of type
$\mathfrak{a}_2$ and the ones whose fix a subalgebra of type
$\mathfrak{g}_2$, corresponding respectively to those whose
centralizer intersected with $\mathcal G_0$ is isomorphic to
$\Int\mathfrak{a}_2$ and $\aut\mathfrak{g}_2$, that is, conjugated
to $\theta'$ and $\theta$ respectively. Next we study which is the
category of $\varphi^2\theta$, which has order 3 since
$\varphi\in\aut C\subset\mathcal G_0$ commutes with $\theta$.

\begin{co}\label{yek}
The centralizer of $\varphi^2\theta$ in $\mathcal G_0$ is
isomorphic to $\Int(\mathfrak{sl}(3,\K))$. In particular
$\varphi^2\theta$ is conjugated to $\theta'$.
\end{co}
\textbf{Proof.} Take $[t]\in\pgo(C,q)^+=\mathcal{G}_0$ such that
$\varphi^2\theta[t]=[t]\varphi^2\theta$. Then $\varphi^2[t^-]=[t]
\varphi^2$ or $[t^-]=[\varphi t\varphi^{-1}]$. A little more
computation reveals that $[t^+]=[\varphi^{-1} t\varphi]$. Thus
there is an admissible triple $(t,\varphi t\varphi^{-1},
\varphi^{-1}t\varphi)$ for $\bar C$. By definition, this means
that $\mu(t)^{-1}t(\bar x\bar y)= \ov{\varphi
t\varphi^{-1}(x)}\,\, \ov{\varphi^{-1}t\varphi(y)}=
t(\varphi^{-1}(x))\ast t(\varphi(y))$ for any $x,y\in C$, if $*$
denotes now the product of $C_{\varphi,\varphi^{-1}}$. Making
$a=\varphi^{-1}(x)$, $b=\varphi(y)$, we get
$\mu(t)^{-1}t(\varphi(\bar a)\varphi^{-1}(\bar b))=t(a)\ast t(b)$
or equivalently $\mu(t)^{-1}t(a\ast b)=t(a)\ast t(b)$ so that some
nonzero multiple of $t$ is an automorphism of $\mathbb P$. This
implies that $C_{\mathcal G_0}(\varphi^2\theta)\subset\aut(\mathbb
P)\cong\Int(\mathfrak{sl}(3,\K))$, so that $C_{\mathcal
G_0}(\varphi^2\theta)$ can not be isomorphic to $G_2$. $\square$

Therefore  the centralizer $C_{\mathcal G}(\theta)$ contains a conjugated of
$\theta'$.

\subsection{On certain quasitori of $\aut(\d4)$ and their induced gradings.}

For any of the previously considered symmetric composition
algebras $S$ (para-Hurwitz or Okubo algebras) we had $\aut(S,*)$
embedded in ${\mathcal{G}}=\aut(\d4)$. If we denote by $f\mapsto
f^\diamond$ the mentioned embedding, we can use it to construct
maximal quasitori in ${\mathcal{G}}$ by  {\sl mixing}  maximal quasitori in
$\aut(S,*)$ together with $\theta_S$, according to Corollary~\ref{key}.

First, consider the case $S=\bar C$ the para-Hurwitz algebra.
Since $\aut(C)=\aut(\bar C)=G_2$, there is a lot of available
information on this group. Let $$
B=\{e_1,e_2,u_1,u_2,u_3,v_1,v_2,v_3\}
$$
  be   the {\em standard basis} of the Cayley algebra $C$,
 defined by
$$\begin{array}{rll}
e_1u_j&=u_j&=u_je_2,\\
e_2v_j&=v_j&=v_je_1,\end{array}\
\begin{array}{rll}
u_iu_j&=v_k&=-u_ju_i ,\\
-v_iv_j&=u_k&=v_jv_i,\end{array}\ \begin{array}{rl} u_iv_i&=e_1,\\
v_iu_i&=e_2,\end{array}
$$
where $e_1$ and $e_2$ are orthogonal idempotents,  $(i,j,k)$ is
any cyclic permutation of $(1,2,3)$, and the remaining relations
are null. Denote by $t_{\a,\b}$ the automorphism of $C$ whose
matrix in the standard basis is the diagonal matrix
$$
\diag(1,1,\alpha,\beta,(\alpha\beta)^{-1},\alpha^{-1},
\beta^{-1},\alpha\beta),
$$
where $\alpha,\beta\in \K^\times$. The set of these automorphisms
is a maximal torus of $\aut(C)$. In particular,
$$
P_1:=\span{\theta,t_{\a,\b}^\diamond\colon \alpha,\beta\in
\K^\times} \cong (\K^\times)^2\times\Z_3
$$
is a maximal quasitorus of ${\mathcal{G}}$. Consider the
$\Z_3$-grading induced by $\theta$ on $\L:=\d4$. This is
$\L=\L_{\bar 0}\oplus\L_{\bar 1}\oplus\L_{\bar 2}$ where the
$\L_{\bar0}$-modules $\L_{\bar1}$ and  $\L_{\bar2}$ in the grading
are dual for the Killing form, hence seven-dimensional (the
fixed subalgebra of $\d4$ by $\theta$ is Der$\,C=\g2=\L_{\bar 0}$.)
Besides they are irreducible (see \cite[Prop 8.6, Ch~8,
p.~138]{Kac}), so that they are isomorphic to the natural
$\g2$-module $C_0$. Since all the components of the $\Z^2$-grading
on $C_0$ are one-dimensional (with generators
$e_1-e_2,u_1,u_2,u_3,v_1,v_2$ and $v_3$), and the same happens for
the non-zero homogeneous components on $\der(C)$ (the root
spaces), the $\Z^2\times\Z_3$-grading on $\L$ induced by $P_1$ is
of type $(12,1)+2(7,0)=(26,1)$.

Now, consider the maximal quasitorus of $\aut(C)$ generated by
 $\{t_{1,-1},t_{-1,1},f\}$, where $f$ is the
automorphism given by the following matrix relative to the
standard basis
$$
f=\begin{pmatrix}
 0 & 1 & 0 & 0 & 0 & 0 & 0 & 0\cr
 1 & 0 & 0 & 0 & 0 & 0 & 0 & 0\cr
 0 & 0 & 0 & 0 & 0 & 1 & 0 & 0\cr
 0 & 0 & 0 & 0 & 0 & 0 & 1 & 0\cr
 0 & 0 & 0 & 0 & 0 & 0 & 0 & -1\cr
 0 & 0 & 1 & 0 & 0 & 0 & 0 & 0\cr
 0 & 0 & 0 & 1 & 0 & 0 & 0 & 0\cr
 0 & 0 & 0 & 0 & -1 & 0 & 0 & 0
\end{pmatrix}.
$$
Again  all the components of the  grading induced on $C_0$ are
one-dimensional,  with generators $e_1-e_2,u_1\pm  v_1,u_2\pm v_2$
and $ u_3\pm v_3$, but all the non-zero homogeneous components of
$\L_{\bar0}$ are Cartan subalgebras of $\g2$, hence
two-dimensional. Thus, the maximal quasitorus of ${\mathcal{G}}$
given by
$$
P_2:=\span{\theta,t_{1,-1}^\diamond,t_{-1,1}^\diamond,f^\diamond }
\cong \Z_2^3\times\Z_3\cong\Z_2^2\times\Z_6
$$
  induces a
$\Z_2^2\times\Z_6$-grading on $\L$ of type
$(0,7)+2(7,0)=(14,7)$.\smallskip

Second, consider the Okubo algebra which is isomorphic to the
pseudo-octonion algebra $\mathbb P$ as previously mentioned.  If
$p\in\GL(3,\K)$, the map
$$
\In(p)\colon \mathbb P\to \mathbb P,\quad \In(p)(x)=pxp^{-1}
$$
is an automorphism of $\mathbb P$.  Take
$$
p_1=  \begin{pmatrix}0 & 0 & 1 \cr
  1 & 0 & 0 \cr
  0 & 1 & 0
  \end{pmatrix}           ,\qquad p_2=\begin{pmatrix}1 & 0 & 0 \cr
  0 & \omega & 0 \cr
  0 & 0 & \omega^2
  \end{pmatrix},
$$
order-three invertible matrices verifying $p_1p_2=\omega p_2p_1$.
Note that
$$\In(p_1)\In(p_2)(x)=p_1p_2x(p_1p_2)^{-1}=\omega
p_2p_1x\omega^{-1}(p_2p_1)^{-1}=\In(p_2)\In(p_1)(x),$$ that is,
$\span{\In(p_1),\In(p_2)}\le \aut(\mathbb P)$ is an abelian
subgroup of automorphisms isomorphic to $\Z_3^2$. Hence
$$
P_3:=\span{\theta',\In(p_1)^\diamond,\In(p_2)^\diamond}\cong\Z_3^3
$$
is another maximal quasitorus of ${\mathcal{G}}$. Consider again
the   $\Z_3$-grading $\L=\L_{\bar 0} \oplus \L_{\bar
1}\oplus\L_{\bar 2}$ induced in $\d4$ by $\theta'$. Then
$\L_{\bar0}$ is $8$-dimensional (isomorphic to
$\mathfrak{sl}(3,\K)$), and $\L_{\bar 1}$ and $\L_{\bar 2}$ are
$\L_{\bar 0}$-dual irreducible modules (again \cite[Prop.~8.6,
Ch.~8, p.~138]{Kac}), hence of dimension 10. More precisely, if
$V$ is a tridimensional vector space, the Lie algebra $\L_{\bar
0}$ is isomorphic to $\mathfrak{sl}(V)$ and the other components
are isomorphic to $\mathfrak{sl}(V)$-modules of type $S^3(V)$ and
$S^3(V^*)$. When considering the $\Z_3$-grading on $V$ given by an
arbitrary basis $\{v_0,v_1,v_2\}$, the map $f_1\colon V\to V$
given by $f_1(v_i)=v_{i+1}$ is extended to $\L$ as
$\In(p_1)^\diamond$, and splits $S^3(V)$ into \lq\lq pieces" of
size $4,3$ and $3$ respectively. The same happens to
$\In(p_2)^\diamond$, which is the extension to $\L$ of $f_2\colon
V\to V$ given by $f_2(v_i)=\omega^i v_{i}$. Both automorphisms
together split $S^3(V)$ in one subspace of dimension $2$ and the
remaining ones of dimension $1$. Thus, the $\Z_3^3$-grading
induced on $\d4$ is of type $(8,0)+2(8,1)=(24,2)$.\smallskip

To summarize, as $\span{ t_{1,-1},t_{-1,1},f } \cong \Z_2^3$ and
$\span{ t_{\a,\b} \colon \alpha,\beta\in \K^\times} \cong
(\K^\times)^2$ are MAD-groups of $\aut\g2$ (the only MAD's, according
to \cite{g2}), and $\span{ \In(p_1),\In(p_2)}\cong\Z_3^2$ is a MAD
of the group of automorphisms of the pseudo-octonions algebra, by
Corollary~\ref{key} we have   that the quasitori $P_i$ are maximal
(and because of their types,  $P_1\cong Q_{13}$, $P_2\cong Q_{12}$ and
$P_3\cong Q_{14}$.)

In order to prove that the $P_i$'s are all the maximal quasitori
(which intersect the $\theta$-component), we prove first that
any other MAD would have some   order 3 outer automorphism.

Consider the maximal torus $T$ of $\mathcal G_0$ induced by the
elements $t_{\lambda, \alpha,\beta,\gamma}\in\go(C,q)$ whose
matrix in the standard basis of $C$ is
\begin{equation}\label{untoro}
\hbox{diag}(\lambda,\lambda^{-1},
\alpha,\beta,\gamma,\alpha^{-1},\beta^{-1},\gamma^{-1}.)
\end{equation}
Straightforward computations reveal that for $t\in T$ one also has
$\bar t,t^+,t^-\in T$ and if $t=t_{\lambda,\alpha,\beta,\gamma}$, then
$\bar t=t_{{\bar\lambda}^{-1},\bar\alpha,\bar\beta,\bar\gamma}$ and $t^\pm= t_{\lambda^\pm,\alpha^\pm,\beta^\pm,\gamma^\pm}$ where

\begin{eqnarray}\label{tmasmenos}
\lambda^-=\frac{1}{\sqrt{\alpha\beta\gamma\lambda}},\quad\alpha^-=
\sqrt{\frac{\lambda
\alpha}{\beta\gamma}},\quad\beta^-=\sqrt{\frac{\lambda
\beta}{\alpha\gamma}},\quad\gamma^-= \sqrt{\frac{\lambda
\gamma}{\alpha\beta}},\cr
\lambda^+=\sqrt{\frac{\alpha\beta\gamma}{\lambda}},\quad\alpha^+=
\sqrt{\frac{\alpha}{\beta\gamma\lambda}},\quad\beta^+=\sqrt{\frac{\beta}
{\alpha\gamma\lambda}},\quad\gamma^+=\sqrt{\frac{\gamma}{\alpha\beta\lambda}}.
\end{eqnarray}

\noindent Notice that $\sigma$, $\theta$ and $\theta'$ belong to
$N(T)$, according to their actions on $T\subset\mathcal{G}_0$
given in (\ref{actionS3}).

\begin{lm}\label{newby}
%\item The subgroup $G^{(\theta)}\cap T^{(\theta)}$ is finite and
%$T=G^{(\theta)}T^{(\theta)}$.
For any $t\in T$ there is some $s\in T$ such that $ts(s^-)^{-1}\in
T^{(\theta)}:=T\cap C_{\mathcal G}(\theta)$.
\end{lm}
Proof. We must take into account that
$t_{\lambda,\alpha,\beta,\gamma}\in T^{(\theta)}$ if and only if
$\lambda=1$ and $\alpha\beta\gamma=1$.  Thus, making
$t=t_{\lambda,\alpha,\beta,\gamma}$ and $s=t_{x,y,u,v}$ the fact
that $ts(s^-)^{-1}\in T^{(\theta)}$ is equivalent to proving that
the equations
$$\lambda x=x^-,\quad \alpha\beta\gamma yuv=u^-v^-y^-,$$ have some solution in
$x$, $y$, $u$ and $v$. But writing $x^-,y^-,u^-$ and $v^-$ as
functions of $x,y,u$ and $v$ according to the relations in
(\ref{tmasmenos}), the resulting equations are
$\alpha\beta\gamma(uvy)^{3/2}=x\sqrt{x}$, $\lambda x
\sqrt{x}=(uvy)^{-1/2}$ which can be solved in $x,y,u,v$ for any
given $\lambda,\alpha,\beta,\gamma$.  $\square$

\begin{pr}\label{hayunodeorden3}
Let $Q$ be a MAD-group of $\mathcal G$ such that
$Q\cap\mathcal{G}_0\theta\ne\emptyset$. Then there is an order 3
outer automorphism in $Q$.
\end{pr}
\textbf{Proof.} Recall that $\mathcal G=\mathcal G_0\cup \mathcal
G_0\sigma\cup \mathcal G_0{\sigma\theta}\cup \mathcal
G_0{\sigma\theta^2}\cup \mathcal G_0\theta\cup \mathcal
G_0{\theta^2}$. If $f\in\mathcal G_0\theta$ has order $3^k\cdot m$
with gcd$(m,3)=1$, then $f^m$ has order $3^k$ and belongs to
$\mathcal G_0\theta$ or to $\mathcal G_0{\theta^2}$, which is
conjugated to $\mathcal G_0\theta$.   Thus we can consider from
the beginning an element
$f\in\mathcal G_0\theta$ with minimum order $3^k$. If $k=1$ we are done, so suppose $k>1$.
%By \cite[Theorem 3.12 (4), p. 115]{enci},
%$t$ can be taken in $T^{(\theta)}$ (conjugating if necessary.)
We must note that there is a maximal toral subgroup $B$ of $Q$ and
a maximal torus $T$ of $\mathcal G$ such that $B\subset T$ and
$Q\subset N(T)$ (see   \cite[Sect.~5]{f4}.) Then $B=T\cap Q$ and
for each $f\in Q\setminus T$ the subquasitorus generated by
$\{f\}\cup(Q\cap T)$ is nontoral. As $T$ is conjugated to the
torus given by  (\ref{untoro}), replace $\theta$, $\theta'$ and
$\sigma$ by their corresponding conjugated automorphisms inside
$N(T)$. Denote by $\pi\colon N(T)\to
N(T)/T\cong\aut\Phi=\mathcal{W}\rtimes S_3$ the canonical
projection ($S_3$ is generated by $\pi(\theta)$ and
$\pi(\sigma)$.) Consider also $\tors{\sigma}:=\{t\in T\mid
\sigma\cdot t=t\}$, for $\cdot$ the action of the Weyl group
$\mathcal W$ on the torus. Notice that $Q\cap
T=\cap_{\sigma\in\pi(Q)}\tors{\sigma}$, by the maximality of $Q$.
Now, $\pi(f)\ne 1$ since $f\notin T$, hence $\pi(f)$ is an order 3
element (since the Weyl group of $\d4$ has no elements of order
$3^n$ for $n>1$.) In particular $\pi(f)$ is conjugated to
$\pi(\theta)$ or $\pi(\theta')$. More concretely, it is conjugated
to $\pi(\theta)$, because every  element in $N(T)$ projecting in
$\pi(\theta')$ has order just 3, following \cite[Lemma~1,
p.~26]{f4}, so in such case $f$ would also have order 3. Thus, by
conjugating by an element in $\mathcal{W}$ if necessary (so we do
not change $T$ nor $N(T)$) we may suppose that
$\pi(f)=\pi(\theta)$, and $f=t\theta\in Q$ for some $t\in T$. Now
conjugating by a suitable element $s\in T$ we have $st\theta
s^{-1}=st(s^-)^{-1}\theta$ and by Lemma~\ref{newby} the toral
element $s$ can be chosen so that $st(s^-)^{-1}\in T^{(\theta)}$.
Thus from here on, we suppose $f=t\theta$ with $t\in
T^{(\theta)}$.

The quasitorus $Q$ being an abelian group, has no element of the
connected components $\mathcal G_0\sigma$, $\mathcal
G_0\sigma\theta$ or $\mathcal G_0\sigma\theta^2$. Therefore
$\pi(Q)=(\pi(Q)\cap \mathcal W)\langle\pi(\theta)\rangle$. Note
that there is  not $g\in Q$ with $\pi(g)\in\mathcal{W}$ of order
3. Because $B\cup\{g\}$ would be non toral, but it is contained in
$T^{(g)}\cup\{g\}$, toral according to Table 2 (and immediate by
Lemma~\ref{toroporalgotoral}.) Therefore
 any element in
$\pi(Q)\cap\mathcal W$ has order a power of 2. If
$\pi(Q)\cap\mathcal W=\{1\}$ we have $\pi(Q)=\span{\pi(\theta)}$
and $Q\cap T=T^{\langle\pi(\theta)\rangle}=T^{(\theta)}\ni t$ so
that $\theta\in Q$ would be an order 3 element.
Otherwise, there is $g\in Q$ with $\pi(g)\in\mathcal{W}$ of order
2. Thus $B\subset T^{\langle\pi(gf)\rangle}$, but for the order 6
elements $\pi(gf)=\pi(g\theta)\in\aut\Phi$ we know by Table~2
that $T^{\langle\pi(gf)\rangle}$ is either  a one-dimensional
torus, or  isomorphic to $\Z_2^2$,  or $\{\id\}$. The last two
options are not possible because $t^3=f^3\in B$ has order
$3^{k-1}$.   Now the
quasitorus generated by $T^{\langle\pi(gf)\rangle}\cup\{g\}$ is
toral
 by Lemma
\ref{toroporalgotoral}, so that the quasitorus generated by
$B\cup\{g\}$ is also toral, a contradiction with the election of
$B$. $\square$
\medskip

Up to this point we have ``mixed" MAD-groups of $\aut(S)$ (where
$S=\bar C$ or $\mathbb P$) with $\theta_S$ so as to obtain
MAD-groups of $\mathcal G=\aut(\d4)$.  Notice that in
$\Int(\mathfrak{a}_2)$ there is another MAD-group, the maximal
torus, but when mixing with $\theta'$ the obtained quasitori will
turn out to be conjugated to $P_1$, as a consequence of
Corollary~\ref{yek}. Thus we arrive at the main theorem.

\begin{te}\label{elimp}
$P_1$, $P_2$ and $P_3$ are, up to conjugation, the only maximal
quasitori of $ \mathcal G$ not contained in $\mathcal G_0\cup
\mathcal G_0\sigma$.
\end{te}
\textbf{Proof.} Let $S=\bar C$ or $\mathbb P$ and consider a
MAD-group $Q$ of $\mathcal G$ containing some element in the
connected component of the element $\theta$. If $Q$ contains some
element conjugated to $\theta'=\theta_\mathbb P$, we can suppose
that $\theta'$ belongs to $Q$. Hence $Q\subset C_{\mathcal
G}(\theta')={(A_2)_0}\cup {(A_2)_0}\theta' \cup
{(A_2)_0}\theta'^2$ applying Proposition~\ref{centr}. Thus, there
is a quasitorus $Q'$ of ${(A_2)_0}=\aut(\mathbb P)$ such that
$Q=\span{Q'\cup\{1,\theta',\theta'^2\}}$. Moreover, the maximality
of $Q$ implies that of $Q'$ in ${(A_2)_0}$.
 There are
only two MAD-groups of ${(A_2)_0}$ (there are four MAD-groups of
$A_2$, according to \cite{gradA2}, but two of them have outer
automorphisms),   and, when we mix them with $\theta'$ we get two
maximal quasitori: $P_3$ and
$$P_4:=\span{\theta',\In(\alpha E_{1,1}+\alpha^{-1} E_{2,2}+
E_{3,3} )^\diamond,\In(E_{1,1}+\beta^{-1} E_{2,2}+\beta E_{3,3}
)^\diamond\colon \alpha,\beta\in \K^{\times}}$$
$$\cong\Z_3\times(\K^{\times})^2,$$
 obtained by joining
the maximal torus of $A_2$ with $\theta'$.

On the other hand, if $Q$ does not contain any element conjugated
to $\theta'$, by Proposition~\ref{hayunodeorden3} there is an
element conjugated to $\theta= \theta_{\bar C}$   in $Q$. Again
$Q$ is a copy of the direct product of a MAD-group of $\aut\g2$
with $\{1,\theta,\theta^2\}$. Notice that there are 2 MAD's of the
group $G_2$ (see \cite{g2}), which provide just $P_1$ and $P_2$
when mixing with $\theta$. But $Q$ cannot be conjugated to $P_2$
because in $P_2$ there is an automorphism conjugated to $\theta'$,
by Corollary~\ref{yek}.

Therefore we have proved that there are only three MAD's in the
 stated conditions, $P_1$, $P_3$ and $P_4$ (in particular, $P_2$ and $P_4$
are conjugated.)
 $\square$

\end{document}